\documentclass[11pt]{article}
\usepackage{anysize}
\marginsize{2.0cm}{2.0cm}{2.0 cm}{2.0 cm}
\usepackage{setspace}

\usepackage{latexsym,amsfonts,amsmath}
\usepackage{graphicx}
\usepackage{color}
\usepackage{epstopdf}
\usepackage[titletoc]{appendix}

\def\liminf{\mathop{\underline{\lim}}}
\def\limsup{\mathop{\overline{\lim}}}

\newtheorem{condition}{Condition}[section]{\bfseries}{\itshape}

{\bfseries}{\itshape}

\newtheorem{theorem}{Theorem}[section]{\bfseries}{\itshape}

\newtheorem{corollary}{Corollary}[section]{\bfseries}{\itshape}

\newtheorem{proposition}{Proposition}[section]{\bfseries}{\itshape}

{\bfseries}{\itshape}

\newtheorem{lemma}{Lemma}[section]{\bfseries}{\itshape}

\newtheorem{remark}{Remark}[section]{\bfseries}{\itshape}

\newtheorem{definition}{Definition}[section]{\bfseries}{\itshape}

{\bfseries}{\itshape}

\begin{document}

\title{On the nonexplosion and explosion for nonhomogeneous Markov pure jump processes}
\author{Yi
Zhang \thanks{Department of Mathematical Sciences, University of
Liverpool, Liverpool, L69 7ZL, U.K.. E-mail: yi.zhang@liv.ac.uk.}}
\date{}
\maketitle

\par\noindent{\bf Abstract:}
In this paper, we obtain new drift-type conditions for nonexplosion and explosion for nonhomogeneous Markov pure jump processes in Borel state spaces. The conditions are sharp; e.g., the one for nonexplosion is necessary if the state space is in addition locally compact and the $Q$-function satisfies weak Feller-type and local boundedness conditions. We comment on the relations of our conditions with the existing ones in the literature, and demonstrate some possible applications.
\bigskip

\par\noindent {\bf Keywords:} Dynkin's formula. Nonhomogeneous Markov pure jump process. Nonexplosion
\bigskip

\par\noindent
{\bf AMS 2000 subject classification:}  60J75. 90C40

\section{Introduction}
In this paper we obtain conditions for nonexplosion and explosion of nonhomogeneous Markov pure jump processes in Borel state spaces. More precisely, we consider the transition rate or say the $Q$-function (or $Q$-matrix in case the state space is denumerable) to be conservative and stable. According to Feller \cite{Feller:1940}, see also \cite{FeinbergShiryayev:2014}, for each $Q$-function, there is a (minimal) transition function solving the Kolmogorov backward equation, giving rise to a (minimal) nonhomogeneous Markov pure jump process. It is this process that is of our interest here. The relevant definitions will be given in the next section.

Markov pure jump processes provide a useful tool for modeling in queueing systems, telecommunication, inventory, and stochastic particle systems etc. In the context of stochastic particle systems, one often needs consider the process to be in an uncountable state space \cite{Eibeck:2003,Wagner:2005}, whereas in the context of continuous-time Markov decision processes, for the sake of optimality in general one needs deal with a nonstationary strategy inducing a nonhomogeneous Markov pure jump process. A stochastic particle system is said to exist if the underlying Markov pure jump process is nonexplosive, see e.g., \cite{Eibeck:2003}.

The study of conditions for nonexplosion of Markov pure jump processes has been done by many authors for many years but mainly for the homogeneous denumerable case, i.e., when the $Q$-matrix does not depend on the time, and the state space is denumerable; such a Markov pure jump process is called a continuous-time Markov chain. A desirable condition would be sharp and at the same time practically applicable. In this connection, for continuous-time Markov chains, Reuter \cite{Reuter:1957} in 1950s obtained a condition, which is necessary and sufficient for the nonexplosion, and in the meanwhile is exclusively in terms of the $Q$-matrix. That condition can be stated as the absence of nontrivial bounded nonnegative eigenvectors of the $Q$-matrix for any positive eigenvalue. Clearly, in practice, it is not easy to verify the Reuter's condition. In response to this, in 1980s Mufa Chen \cite{Chen:1986} obtained a verifiable condition for nonexplosion in terms of the existence of a drift function with respect to the $Q$-function. The sufficiency of this condition was justified for the homogeneous but general state space case in \cite{Chen:1986,Chen:2004}. Very recently, Spieksma \cite{Spieksma:2015} showed that Chen's condition in \cite{Chen:1986} is actually necessary for the continuous-time Markov chain. In the survey by Mufa Chen \cite{Chen:2015}, the necessity of that condition in the homogeneous but general state space case was said to be unclear, see the proof of Theorem 4 as well as the last paragraph of the text therein. (In the present paper the necessity of Chen's condition is justified if the state space is a locally compact Borel space, and the $Q$-function satisfies some weak-Feller type conditions.) For the other relevant results for the homogeneous case, we refer the interested reader to the book \cite{Anderson:1991} and the recent survey \cite{Chen:2015} in case the state space is denumerable, and the monograph \cite{Chen:2004} in the case of a more general state space. See also \cite{Menshikov:2014} for some more recent progresses, where more specific continuous-time Markov chains were considered.

Less is known for the nonhomogeneous case. In case the state space is denumerable, and in addition the $Q$-function is continuous in time, a drift condition for the nonexplosion was independently obtained by Chow and Khasminskii \cite{Chow:2011} recently, see Theorem 1 therein, and Junli Zheng and Xiaogu Zheng \cite{Zheng:1987,Zheng:1993}. In Meyn and Tweedie \cite{MeynTweedie:1993III}, a sufficient condition for the nonexplosion of the general homogeneous continuous-time Markov process in a locally compact seperable metric state space was obtained. The argument in \cite{MeynTweedie:1993III} and the one in \cite{Chow:2011} are similar; in the latter work, the homogenization technique was used. No discussions were provided in \cite{Chow:2011,MeynTweedie:1993III,Zheng:1987,Zheng:1993} on the sharpness of their conditions.

The main contribution of the present paper is that for nonhomogeneous Markov pure jump processes in Borel state spaces we obtain new and sharp drift-type conditions for nonexplosion, where the sharpness is understood in the sense that our condition is necessary in case the state space is locally compact seperable metrizable and some weak Feller-type and local boundedness conditions are satisfied. As an application of this condition, a generalization of the condition in Chow and Khasminskii \cite{Chow:2011} and Junli Zheng and Xiaogu Zheng \cite{Zheng:1987,Zheng:1993} is obtained for the general Borel state space case. Our arguments are different from those of \cite{Chow:2011,MeynTweedie:1993III,Zheng:1987,Zheng:1993}, which, roughly speaking, are based on the study of a sequence of truncated continuous-time models with the truncation being done on the state space. The arguments here are mainly based on the connection between the nonexplosion of the continuous-time process and the stability of an underlying discrete-time Markov process. In the homogeneous denumerable case, this idea was employed by Spieksma in \cite{Spieksma:2015}. It was noted in \cite{Spieksma:2012} that for a continuous-time Markov chain, whether Dynkin's formula is applicable to a drift function depends on whether a transformed continuous-time Markov chain is nonexplosive. We show that the similar result holds for the nonhomogeneous Markov pure jump process in the Borel space. This thus demonstrates an application and the practical importance of a sharp condition for the nonexplosion.

The rest of this paper is organized as follows. Basic notations and definitions about nonhomogeneous Markov pure jump processes are summarized in Section \ref{ZY2015:SecDef}. Our conditions for nonexplosion and explosion along with their relation to the existing ones are presented in Section \ref{ZY2015:SecNonExp}. In Section \ref{ZY2015:SecAppl} we study the relation between the applicability of Dynkin's formula to a drift function and the nonexplosion of a transformed process, and its consequences. We finish the paper with a conclusion in Section \ref{ZY2015ConcSec}. The auxiliary statements and proofs are postponed to the appendix.

\section{Definitions, notations and basic facts}\label{ZY2015:SecDef}

In what follows, let $S$ be a nonempty Borel state space endowed with its Borel $\sigma$-algebra ${\cal B}(S).$ The term ``measurability'' is always understood in the Borel sense. The Dirac measure concentrated on the singleton $\{x\}$ is denoted as $\delta_{x}(dy)$.

As in Feinberg, Mandava and Shiryaev \cite{FeinbergShiryayev:2014}, we adopt the following definition throughout this paper.
\begin{definition}
A (Borel-measurable) signed kernel $q(dy|x,s)$ on ${\cal B}(S)$ from $S\times [0,\infty)$ is called a (conservative stable) $Q$-function on the Borel space $S$ if the following conditions are satisfied.
\begin{itemize}
\item[(a)] For each $s\ge 0$, $x\in S$ and $\Gamma\in {\cal B}(S)$ with $x\notin \Gamma,$
\begin{eqnarray*}
\infty>q(\Gamma|x,s)\ge 0.
\end{eqnarray*}
\item[(b)]  For each $(x,s)\in S\times [0,\infty),$
\begin{eqnarray*}
q(S|x,s)=0.
\end{eqnarray*}
\item[(c)] For each $x\in S,$
\begin{eqnarray*}
\sup_{s\ge 0}\left\{q(S\setminus \{x\}|x,s)\right\}<\infty.
\end{eqnarray*}
\end{itemize}
\end{definition}

Given a $Q$-function $q$ on $S$, we introduce the kernel $\tilde{q}(dy|x,s)$  on ${\cal B}(S)$ from $S\times [0,\infty)$ defined by
\begin{eqnarray*}
\tilde{q}(\Gamma|x,s):=q(\Gamma\setminus\{x\}|x,s)
\end{eqnarray*}
for each $x\in S$, $s\in [0,\infty)$ and $\Gamma\in {\cal B}(S)$, and the function
\begin{eqnarray*}
q_x(s):=-q(\{x\}|x,s),~\forall~x\in S,~s\in [0,\infty).
\end{eqnarray*}

For each (conservative stable) $Q$-function $q$ on a Borel space $S$, it is convenient to agree on that
\begin{eqnarray}\label{ZY2015:58}
q_x(\infty):=0=:\tilde{q}(\Gamma|x,\infty)
\end{eqnarray}
for each $x\in S$ and $\Gamma\in {\cal B}(S).$

For each $\Gamma\in {\cal B}(S)$, $x\in S$, $s,t\in [0,\infty)$ and $s\le t,$ we define
\begin{eqnarray}\label{ZY2015:07}
&&P_q^{(0)}(s,x,t,\Gamma):= \delta_{x}(\Gamma) e^{-\int_s^t q_x(v)dv},\nonumber\\
&&P_q^{(n+1)}(s,x,t,\Gamma):= \int_s^t e^{-\int_s^u q_x(v)dv} \left(\int_S  P_q^{(n)}(u,z,t,\Gamma)\tilde{q}(dz|x,u)\right)du,~\forall~ n=0,1,\dots.
\end{eqnarray}

It is well known or otherwise easy to check that the following assertions hold;
\begin{itemize}
\item[(a)] for each $s,t\in [0,\infty),$ $s\le t$ and $x\in S$, $P_q^{(n)}(s,x,t,dy)$ is a $[0,1]$-valued measure on ${\cal B}(S)$;
\item[(b)] for each $\Gamma \in {\cal B}(S)$, $(s,x,t)\rightarrow P_q^{(n)}(s,x,t,\Gamma)$ is jointly measurable in $x\in S$ and $s,t\in[0,\infty),$ $s\le t.$
\end{itemize}

Define
\begin{eqnarray}\label{ZY2015:01}
P_q(s,x,t,\Gamma):=\sum_{n=0}^\infty P_q^{(n)}(s,x,t,\Gamma)
\end{eqnarray}
for each $x\in S$, $s,t\in [0,\infty),$ $s\le t$, and $\Gamma\in {\cal B}(S)$.

The following result is due to Feller \cite{Feller:1940}; see also \cite{FeinbergShiryayev:2014,YeGuo:2015}.
\begin{proposition}\label{ZY2015:Prop02}
\par\noindent(a)
For each $s,t\in [0,\infty),$ $s\le t$ and $x\in S$, $P_q (s,x,t,dy)$ is a $[0,1]$-valued measure on ${\cal B}(S)$.
\par\noindent(b) For each $\Gamma \in {\cal B}(S)$, $(s,x,t)\rightarrow P_q(s,x,t,\Gamma)$ is measurable (jointly) in $x\in S$ and $s,t\in[0,\infty),$ $s\le t.$

(Thus, from (a) and (b), $P_q(s,x,t,dy)$ is a transition function.)
\par\noindent (c) The transition function $P_q(s,x,t,dy)$ satisfies the following Chapman-Kolmogorov equation;
\begin{eqnarray*}
P_q(s,x,t,\Gamma)=\int_S P_q(s,x,u,dy)P_q(u,y,t,\Gamma),~\forall~ x\in S,~ s,u,t\in[0,\infty),~ s<u<t,~\Gamma\in {\cal B}(S).
\end{eqnarray*}
\end{proposition}

Throughout this paper, the transition function $P_q(s,x,t,dy)$ defined by (\ref{ZY2015:01}) is called the $q$-transition function (to signify its dependence on the $Q$-function $q$). In the literature, it is also called the Feller's minimal $q$-transition function, where the minimality is in the sense that it is the minimal nonnegative solution to the Kolmogorov backward equation; see Theorem 3.2 of \cite{FeinbergShiryayev:2014}. The latter fact is not needed in this paper, whereas the following Kolmogorov forward equation will be used.

\begin{proposition}\label{ZY2015:Prop03}
For each $x\in S,$ $s,t\in[0,\infty)$, $s\le t$ and $\Gamma\in {\cal B}(S)$ satisfying
\begin{eqnarray}\label{ZY2015:28}
\sup_{x\in \Gamma,~s\ge 0}q_x(s)<\infty,
\end{eqnarray}
the following holds;
\begin{eqnarray}\label{ZY2015:63}
&&P_q(s,x,t,\Gamma)=\delta_x(\Gamma)+ \int_s^t \left(\int_S \tilde{q}(\Gamma|y,u)P_q(s,x,u,dy)\right)du -\int_s^t \left(\int_\Gamma q_y(u)P_q(s,x,u,dy)\right)du \nonumber\\
&=&\delta_x(\Gamma)+ \int_s^t \left(\int_S q(\Gamma|y,u)P_q(s,x,u,dy)\right)du.
\end{eqnarray}
\end{proposition}
\par\noindent\textit{Proof.} See Corollary 4.2 of \cite{FeinbergShiryayev:2014}. $\hfill\Box$\bigskip

It is worthwhile to point out that in general (\ref{ZY2015:63}) does not hold for all $\Gamma\in {\cal B}(S)$; the statements of Lemma 1(ii) and Theorem 2(ii) of \cite{YeGuo:2015} are not quite accurate.

For the subsequent discussions, the connection between the $q$-transition function $P_q(s,x,t,dy)$ and a (nonhomogeneous) pure jump Markov process is briefly described as follows. (The term of nonhomogeneous is traditional, and it is understood as not necessarily homogeneous.)

Let $x_\infty\notin S$ be some fixed isolated point. Denote $S_\infty:=S\bigcup\{x_\infty\}$.  Let the sample space $\Omega$ be the collection of sequences $(x_0,\theta_1,x_1,\theta_2,\dots)\in (S_\infty\times (0,\infty])^\infty$ such that
\begin{itemize}
\item[(a)] $x_0\in S$, and
\item[(b)] $\theta_n=\infty$ for some $n=1,2,\dots$, if and only if $x_m=x_\infty$ and $\theta_{m+1}=\infty$ for each $m=n,n+1,\dots.$
\end{itemize}
The $\sigma$-algebra ${\cal F}$ on $\Omega$ is the trace of the product $\sigma$-algebra of $(S_\infty\times (0,\infty])^\infty$. The generic notation $\omega\in \Omega$ is in use but typically omitted;
so that in what follows, without special reference the coordinate variables are often written as
\begin{eqnarray*}
x_n(\omega)=x_n,~ \theta_{n+1}(\omega)=\theta_{n+1}.
\end{eqnarray*}
(The abuse of notation here is conventional.)
We define
\begin{eqnarray*}
t_n:= \sum_{i=1}^n \theta_i,~\forall~n=1,2,\dots,~t_0:=0.
\end{eqnarray*}
Below $\theta_n$ and $t_n$ will be interpreted as the $n$th sojourn time and jump moment. The accumulation point of the monotone nondecreasing sequence $(t_n)$ is denoted by $t_\infty.$

The collection $(t_n,x_n)$ is a marked point process. Let the filtration $\{{\cal F}_t\}$ be the internal history of $(t_n,x_n)$. It is known, see e.g., \cite{Jacod:1975} or Chapter 4 of \cite{Kitaev:1995}, that for each fixed initial distribution $\gamma$ on ${\cal B}(S)$, there is a unique probability $\mathbb{P}^q_\gamma$ on $(\Omega, {\cal F})$ such that the following assertions hold.
\begin{itemize}
\item[(a)] $\mathbb{P}^q_\gamma(x_0(\omega)\in \Gamma)=\gamma(\Gamma)$ for each $\Gamma\in {\cal B}(S)$.
\item[(b)] The random measure on ${\cal B}((0,\infty)\times S)$ defined by
\begin{eqnarray*}
\nu(\omega,~dt\times dx):= \sum_{n=0}^\infty I\{t_n<t\le t_{n+1}\}\tilde{q}(dx|x_n,t)dt
\end{eqnarray*}
is the dual predictable projection of the random measure of the marked point process $(t_n,x_n)$.
\end{itemize}
If the initial distribution $\gamma$ is concentrated on a singleton say $\{x\}\subseteq S,$ then $\mathbb{P}_\gamma^q$ is written as $\mathbb{P}_x^q.$ The expectation with respect to $\mathbb{P}_\gamma^q$ and $\mathbb{P}_x^q$ are denoted by $\mathbb{E}_\gamma^q$ and $\mathbb{E}_x^q.$ When it is not important to emphasize on the initial distribution, we simply write $\mathbb{P}^q$ and $\mathbb{E}^q.$

Now we define the process of interest $(X_t,~t\ge 0)$ on $(\Omega, {\cal F})$ with the filtration $\{{\cal F}_t\}$ by
\begin{eqnarray}\label{ZY2015:50}
&&X_t(\omega)=\left\{
\begin{array}{ll}
      x_n, & \mbox{if~}  t_n\le t<t_{n+1}, \\
     x_\infty, & \mbox{if~} t_\infty \le t,
\end{array}
\right.
\end{eqnarray}
 for each $\omega\in \Omega.$

\begin{proposition}\label{ZY2015:Prop01}
For each $u,t\in (0,\infty),$ $u\le t$ and $\Gamma\in{\cal B}(S)$,
\begin{eqnarray*}
\mathbb{P}^q(X_t\in \Gamma|{\cal F}_u)=\mathbb{P}^q(X_t\in \Gamma|X_u)=P_q(u,X_u,t,\Gamma)
\end{eqnarray*}
almost surely with respect to $\mathbb{P}^q$ on $\{u<t_\infty\}.$
\end{proposition}
\par\noindent\textit{Proof.} See Theorem 2.2 of \cite{FeinbergShiryayev:2014}. $\hfill\Box$\bigskip

Thus, $(X_t)$ is a (nonhomogeneous) pure jump Markov process with the transition function $P_q.$

\begin{definition}
Consider the $Q$-function $q$ on $S$. The (nonhomogeneous) pure jump process $(X_t)$ with the transition function $P_q$ is called nonexplosive if
\begin{eqnarray*}
\mathbb{P}_x^q(t_\infty=\infty)=1,~\forall~x\in S.
\end{eqnarray*}
If there exists some $x\in S$ such that $\mathbb{P}_x^q(t_\infty=\infty)<1$, then the process $(X_t)$ is called explosive.
\end{definition}

One objective of this paper is to obtain sharp conditions for nonhomogeneous Markov pure jump processes. To this end, we firstly present the next simple but useful observation concerning the equivalent definitions of nonexplosion.
\begin{lemma}\label{ZY2015:Lemma07}
Consider the $Q$-function $q$ on $S$. Then the following assertions are equivalent.
\begin{itemize}
\item[(a)] The (nonhomogeneous) pure jump process $(X_t)$ with the transition function $P_q$ is nonexplosive.
\item[(b)] For each $x\in S$ and $t\ge 0$, $P_q(0,x,t,S)=1.$
\item[(c)] For each $x\in S,$ and almost all $t\ge0$, $P_q(0,x,t,S)=1.$
 \item[(d)] For each $x\in S$ and $0\le s\le t$, $P_q(s,x,t,S)=1.$
 \end{itemize}
\end{lemma}
 \par\noindent\textit{Proof.} See the appendix. $\hfill\Box$\bigskip

To finish this section, for a $Q$-function $q$ on (the Borel space) $S$ and a constant $c\in[0,\infty),$ let us introduce the $f$-transformed $Q$-function as follows, where $f$ is a $(0,\infty)$-valued measurable function on $S$ satisfying
\begin{eqnarray*}
\int_S f(y)\tilde{q}(dy|x,s)<\infty,~\int_S f(y)q(dy|x,s)\le c f(x)~\forall~x \in S, ~s\ge 0.
\end{eqnarray*}
Let \begin{eqnarray*}
S_\delta:=S\bigcup \{\delta\}
\end{eqnarray*}
with $\delta\notin S$ being an isolated point and satisfies $\delta\ne x_\infty$. (Here we choose a different isolated point because we want to reserve $x_\infty$ as the isolated point in the definition of $(X_t)$, see (\ref{ZY2015:50}), when the state space is $S_\delta.$) We define the (conservative and stable) $Q$-function $q^f$ on $S_\delta$ (instead of $S$) by
\begin{eqnarray}\label{ZY2015:08}
&&q^f(\Gamma|x,s)=\left\{
\begin{array}{ll}
      \frac{\int_{\Gamma}f(y)q(dy|x,s)}{f(x)}, & \mbox{if~}  x\in S,~ \Gamma\in {\cal B}(S),~x\notin \Gamma; \\
     c-\frac{\int_{S}f(y)q(dy|x,s)}{f(x)}, & \mbox{if~} x\in S,~ \Gamma=\{\delta\};\\
     0, & \mbox{if~} x=\delta,~ \Gamma=S.
\end{array}
\right.
\end{eqnarray}
for each $s\ge 0.$  Therefore, the isolated state $\delta$ is absorbing, and the transition rate at $x\in S$ is given by
\begin{eqnarray*}
q^f_x(s)=c+q_x(s),~\forall~s\ge 0.
\end{eqnarray*}
Below we reserve the isolated point $\delta$ only for the $f$-transformed model.
The $q^f$-transition function on $S_\delta$ is denoted by $P_{q^f}(s,x,t,dy)$. The nonhomogeneous pure jump Markov process with the transition function $P_{q^f}$ is still denoted by $(X_t)$ despite its state space is different now, but its corresponding probability and expectation will be signified with the superscript $q^f.$

The notion of the $f$-transformed $Q$-function for a continuous-time Markov chain appeared in Anderson \cite{Anderson:1991} and Spieksma \cite{Spieksma:2012}. The relationship between the $q$-transition function and the $q^f$-transition function is given by Lemma \ref{ZY2015:09} in the appendix.

\section{Explosion and nonexplosion}\label{ZY2015:SecNonExp}

Consider a $Q$-function $q$ on $S$. The objective of this section is to provide conditions for the nonexplosion and explosion of the nonhomogeneous Markov pure jump process $(X_t)$ in a Borel state space $S$ with the transition function $P_q$. The simpler homogeneous case is commented at the end of this section.

It follows from Lemma \ref{ZY2015:Lemma07} that the process $(X_t)$ with the transition function $P_q$ is nonexplosive if and only if for some constant $\alpha\in(0,\infty),$
\begin{eqnarray}\label{ZY2015:55}
\alpha\int_0^\infty e^{-\alpha t}P_q(v,x,v+t,S)dt=1,~\forall~x\in S,~v\in[0,\infty).
\end{eqnarray}

Let $\Delta\notin S$ be an isolated point such that $\Delta\ne \delta$ and $\Delta\ne x_\infty.$ Now consider the (homogeneous) discrete-time Markov chain $(\tilde{t}_n,\tilde{x}_n)$ in the state space $([0,\infty)\times S)\bigcup(\{\infty\}\times\{\Delta\})$ with the transition probability given by
\begin{eqnarray}\label{ZY2015:53}
&&p^\alpha(\Gamma_1\times \Gamma_2|v,x)\nonumber\\
&:=&\left\{
\begin{array}{ll}
     \int_{\Gamma_1\bigcap(v,\infty)} e^{-\int_0^{t-v} (\alpha+q_{x}(s+v))ds}  \tilde{q}(\Gamma_2|x,t)dt, & \mbox{if~}  x\in S,v\in [0,\infty), \Gamma_1\in {\cal B}([0,\infty)),~\Gamma_2\in {\cal B}(S); \\
     1- \int_0^\infty e^{-\alpha t -\int_0^t q_{x}(s+v)ds} q_{x }(t+v)dt, & \mbox{if~} x\in S,~v\in [0,\infty),~ \Gamma_1=\{\infty\},~\Gamma_2=\{\Delta\};\\
     1, & \mbox{if~} x=\Delta,~v=\infty,~ \Gamma_1=\{\infty\},~ \Gamma_2=\{\Delta\}.
\end{array}
\right.\nonumber\\
\end{eqnarray}
(The transition probability on other states to other sets can be automatically calculated based on the above definitions.) Clearly, the state $(\infty,\Delta)$ is an isolated (cemetery) point. Let $\tilde{P}_{(v,x)}$ be the corresponding probability measure for this discrete-time Markov chain with the transition probability $p^\alpha$ and the initial state $(v,x)\in ([0,\infty)\times S)\bigcup (\{\infty\}\times \{\Delta\}).$

One can see that for each $x\in S$ and $v\in [0,\infty),$
\begin{eqnarray}\label{ZY2015:56}
&&\alpha\int_0^\infty e^{-\alpha t}P_q(v,x,v+t,S)dt=\tilde{E}_{(v,x)}\left[\sum_{n=0}^\infty I\{\tilde{x}_{n}\in S\}\int_0^\infty \alpha e^{-\alpha t}e^{-\int_0^t q_{\tilde{x}_n}(s+\tilde{t}_n)ds}  dt\right]\nonumber\\
&=&\tilde{E}_{(v,x)}\left[\sum_{n=0}^\infty I\{\tilde{x}_{n}\in S\}\left(1- \int_0^\infty e^{-\alpha t -\int_0^t q_{\tilde{x}_n}(\tilde{t}_n+s)ds} q_{x_n}(\tilde{t}_n+t)dt\right)\right]\nonumber\\
&=&\tilde{E}_{(v,x)}\left[\sum_{n=0}^{\infty}I\{\tilde{x}_{n}\in S\}p^\alpha(\{\infty\}\times \{\Delta\}|\tilde{t}_n,\tilde{x}_n)\right]\nonumber\\
&=&\sum_{n=0}^{\infty}\tilde{E}_{(v,x)}\left[I\{\tilde{x}_{m}\in S,~\tilde{t}_m\in [0,\infty),~m=0,1,\dots,n\}I\{(\tilde{t}_{n+1},\tilde{x}_{n+1})=(\infty,\Delta) \}\right].
\end{eqnarray}
Despite this argument is standard, to improve the readability, we carefully verify the previous equalities for the case of $v=0$ as follows. For the constant function $x\in S\rightarrow \alpha\in(0,\infty)$, also written as $\alpha,$ consider the $Q$-function $q^\alpha$ on $S_\delta$ as introduced by (\ref{ZY2015:08}), where $c$ is replaced with $\alpha.$
By Lemma \ref{ZY2015:09},
\begin{eqnarray}\label{ZY2015:51}
\alpha\mathbb{E}_x^q\left[\int_0^\infty e^{-\alpha t}I\{X_t\in S\} dt\right]=\alpha\mathbb{E}_x^{q^\alpha}\left[\int_0^\infty I\{X_t\in S\} \right],~\forall~x\in S.
 \end{eqnarray}
Let us rewrite the right hand side of (\ref{ZY2015:51}) as the expected total cost of a homogeneous discrete-time Markov process as follows; for each $x\in S,$
\begin{eqnarray}\label{ZY2015:52}
&&\alpha\mathbb{E}_x^{q^\alpha}\left[\int_0^\infty I\{X_t\in S\} \right]=\alpha \sum_{n=0}^\infty\mathbb{E}_x^{q^\alpha}\left[\int_{t_n}^{t_{n+1}} I\{x_n\in S\}dt \right]\nonumber\\
&=&\alpha\left.\sum_{n=0}^\infty\mathbb{E}_x^{q^\alpha}\left[\mathbb{E}_x^{q^\alpha}\left[\int_{t_n}^{t_{n+1}} I\{x_n\in S\}dt \right|{\cal F}_{t_n}\right]\right]=\sum_{n=0}^\infty\mathbb{E}_x^{q^\alpha}\left[I\{x_{n}\in S\}\int_0^\infty \alpha e^{-\alpha t}e^{-\int_0^t q_{x_n}(s+t_n)ds}  dt\right]\nonumber\\
&=&\sum_{n=0}^\infty\mathbb{E}_x^{q^\alpha}\left[I\{x_{n}\in S\}\left(1- \int_0^\infty e^{-\alpha t -\int_0^t q_{x_n}(t_n+s)ds} q_{x_n}(t+t_n)dt\right)\right]
\end{eqnarray}
where ${\cal F}_{t_n}=\sigma\{(t_m,x_m),~1\le m\le n\}$, c.f. Theorem 4.13 of Kitaev and Rykov \cite{Kitaev:1995}, and the last equality is by integration by parts.

For each fixed $x\in S,$ one can observe that the process $(\tilde{t}_n,\tilde{x}_n)$ under $\tilde{P}_{(0,x)}$ is rather similar to the marked point process $(t_n,x_n)$ under $\mathbb{P}_x^{q^\alpha}$. The differences are that the absorbing (cemetery) state for $x_n$ is now replaced by $\Delta$, and we allow the transition from $[0,\infty)\times S$ to $\{\infty\}\times \{\Delta\}$ for $(\tilde{t}_n,\tilde{x}_n)$ with positive probability at each state, whereas $t_n$ is finite with full probability (given $x_{n-1}\in S$) because of the presence of $\alpha>0$ in the $Q$-function $q^\alpha$.  However, due to the presence of $I\{x_{n}\in S\}$, the last expression of  (\ref{ZY2015:52}) can be rewritten as
\begin{eqnarray*}
&&\sum_{n=0}^\infty\mathbb{E}_x^{q^\alpha}\left[I\{x_{n}\in S\}\left(1- \int_0^\infty e^{-\alpha t -\int_0^t q_{x_n}(t_n+s)ds} q_{x_n}(t_n+t)dt\right)\right]\nonumber\\
&=&\tilde{E}_{(0,x)}\left[\sum_{n=0}^{\infty}I\{\tilde{x}_{n}\in S\}p^\alpha(\{\infty\}\times \{\Delta\}|\tilde{t}_n,\tilde{x}_n)\right]\nonumber\\
&=&\tilde{E}_{(0,x)}\left[\sum_{n=0}^{\infty}I\{\tilde{x}_{m}\in S,~\tilde{t}_m\in [0,\infty),~m=0,1,\dots,n\}I\{(\tilde{t}_{n+1},\tilde{x}_{n+1})=(\infty,\Delta) \}\right]
\end{eqnarray*}
for each $x\in S$. Thus, (\ref{ZY2015:56}) is verified for the case of $v=0.$ The similar argument applies to the general case which only involves extra notations, and is thus omitted.

The relation (\ref{ZY2015:56}) is fundamentally important in this section because it gives two views of the explosion and nonexplosion of the process $(X_t)$. Indeed, the right hand side expressions can be viewed as the expected total cost well studied in Markov decision processes, giving rise to the first view, as well as the probability of getting absorbed at $\{\infty\}\times \{\Delta\}$ of a discrete-time Markov chain welll studied in stability, giving the second. The first view easily gives a necessary and sufficient condition for the explosion, whereas the second view easily leads to a sharp drift condition for nonexplosion. We elaborate on this in the rest of this section.

\begin{theorem}\label{ZY2015Theorem05}
Consider the $Q$-function $q$ on $S$. Then the (nonhomogeneous) Markov pure jump process $(X_t)$ with the transition function $P_q$ is nonexplosive if and only if one of the following two equivalent assertions hold.
 \begin{itemize}
\item[(a)] There exists some constant $\alpha\in(0,\infty)$ such that the (homogeneous) discrete-time Markov chain $(\tilde{t}_n,\tilde{x}_n)$ with the transition probability $p^\alpha$ given by (\ref{ZY2015:53}) is Harris recurrent with respect to the irreducibility measure $\psi$ concentrated on the singleton $\{\infty\}\times \{\Delta\}.$ (See e.g., \cite{MeynTweedie:1993} for the relevant definitions about the stability of discrete-time Markov chains.)

\item[(b)] For each $\alpha\in (0,\infty)$, there is no nontrivial $[0,1]$-valued measurable function $U$ on $[0,\infty)\times S$ satisfying
 \begin{eqnarray}\label{ZY2015:59}
U(v,x)=\int_0^\infty \int_S U(v+t,y)\tilde{q}(dy|x,v+t)e^{-\alpha t -\int_0^t q_{x}(s+v)ds}dt,~\forall~x\in S,~v\in[0,\infty).
 \end{eqnarray}
 Here and below, by saying $U$ is nontrivial, we mean that there exists some $(v,x)$ such that $U(v,x)\ne 0.$
 \end{itemize}
\end{theorem}
\par\noindent\textit{Proof.}
(a) The summand in the last expression of (\ref{ZY2015:56}) is the first entrance probability of the process $(\tilde{t}_n,\tilde{x}_n)$ into $\{\infty\}\times\{\Delta\}$. The statement immediately follows from this observation together with (\ref{ZY2015:55}), (\ref{ZY2015:53}) and (\ref{ZY2015:56}).

(b) Let $\alpha>0$ be arbitrarily fixed. The first expression of (\ref{ZY2015:56}) is the expected total cost for the discrete-time Markov chain $(\tilde{t}_n,\tilde{x}_n)$ with the cost function given by
\begin{eqnarray*}
(v,x)\in ([0,\infty)\times S)\bigcup (\{\infty\}\times \{\Delta\})\rightarrow I\{x\in S\}\int_0^\infty \alpha e^{-\alpha t}e^{-\int_0^t q_{x}(s+ v)ds}dt;
\end{eqnarray*}
recall the convention adopted in (\ref{ZY2015:58}). It is well known, see e.g., Propositions 9.8 and 9.10 of Bertsekas and Shreve \cite{Bertsekas:1978}, that the function
\begin{eqnarray*}
&&(v,x)\in[0,\infty)\times S\rightarrow \tilde{E}_{(v,x)}\left[\sum_{n=0}^\infty I\{\tilde{x}_{n}\in S\}\int_0^\infty \alpha e^{-\alpha t}e^{-\int_0^t q_{\tilde{x}_n}(s+\tilde{t}_n)ds}  dt\right]\\
&=&\alpha\int_0^\infty e^{-\alpha t}P_q(v,x,v+t,S)dt
 \end{eqnarray*}
 (c.f. (\ref{ZY2015:56})) is the minimal solution to the following dynamic programming-type equation
\begin{eqnarray*}
W(v,x)&=& \int_0^\infty \alpha e^{-\alpha t}e^{-\int_0^t q_{x}(s+v)ds}dt+\int_0^\infty\int_S W(t+v,y) \tilde{q}(dy|x,v+t)e^{-\alpha t -\int_0^t q_x(v+s)ds}dt,\\
&&~\forall~x\in S,~v\in[0,\infty),
\end{eqnarray*}
out of the class of $[0,1]$-valued measurable functions on $[0,\infty)\times S.$ It is easy to check that a $[0,1]$-valued measurable function $W$ on $[0,\infty)\times S$ is a solution to the above equation if and only if $U=1-W$ is a solution to (\ref{ZY2015:59}). Consequently, \begin{eqnarray*}
(v,x)\in[0,\infty)\times S\rightarrow 1- \alpha\int_0^\infty e^{-\alpha t}P_q(v,x,v+t,S)dt
 \end{eqnarray*}
 is the maximal solution to (\ref{ZY2015:59}) out of the class of $[0,1]$-valued measurable functions on $[0,\infty)\times S.$ The statement now follows.
$\hfill\Box$
\bigskip

Part (b) of the previous theorem was stated in the Chinese literature, see e.g., Theorem 7.2 of Chapter 3 in \cite{Hu:2013}, where the continuity of the function $t\in[0,\infty)\rightarrow q(\Gamma|x,t)$ for each $x\in S$ and $\Gamma\in {\cal B}(S)$ was assumed. The proof here is essentially the same as the one of \cite{Hu:2013}, but we emphasize more explicitly on the connections with the underlying discrete-time Markov chain. In practice, it is not easy to apply directly the condition in part (b) of Theorem \ref{ZY2015Theorem05} to checking for the nonexplosion. Nevertheless, the negation of its equivalent form provides a more applicable condition for the explosion of the process $(X_t)$. The statement is as follows.
\begin{corollary}\label{ZY2015Corollary2}
Consider the $Q$-function $q$ on $S$. The (nonhomogeneous) Markov pure jump process $(X_t)$ with the transition function $P_q$ is explosive if and only if for some $\alpha>0,$ there is a nontrivial bounded nonnegative measurable function $U$ satisfying the following inequality:
 \begin{eqnarray}\label{ZY2015:62}
 U(v,x)\le\int_0^\infty \int_S U(v+t,y)\tilde{q}(dy|x,v+t)e^{-\alpha t -\int_0^t q_{x}(s+v)ds}dt,~\forall~x\in S,~v\in[0,\infty).
 \end{eqnarray}
\end{corollary}

\par\noindent\textit{Proof.} The statement directly follows from the proof of part (b) of Theorem \ref{ZY2015Theorem05} and known facts in Markov decision processes. In greater detail, by Proposition 9.10 of Bertsekas and Shreve \cite{Bertsekas:1978}, the function
\begin{eqnarray*}
&&(v,x)\in[0,\infty)\times S\rightarrow \tilde{E}_{(v,x)}\left[\sum_{n=0}^\infty I\{\tilde{x}_{n}\in S\}\int_0^\infty \alpha e^{-\alpha t}e^{-\int_0^t q_{\tilde{x}_n}(s+\tilde{t}_n)ds}  dt\right]
 \end{eqnarray*}
is the minimal solution to the following inequality
\begin{eqnarray*}
W(v,x)&\ge& \int_0^\infty \alpha e^{-\alpha t}e^{-\int_0^t q_{x}(s+v)ds}dt+\int_0^\infty\int_S W(t+v,y) \tilde{q}(dy|x,v+t)e^{-\alpha t -\int_0^t q_x(v+s)ds}dt,\\
&&\forall~x\in S,~v\in[0,\infty)
\end{eqnarray*}
out of the class of $[0,1]$-valued measurable functions on $[0,\infty)\times S.$ Now,
as in the proof of  Theorem \ref{ZY2015Theorem05}(b), it follows that the function
\begin{eqnarray*}
(v,x)\in[0,\infty)\times S\rightarrow 1- \alpha\int_0^\infty e^{-\alpha t}P_q(v,x,v+t,S)dt
 \end{eqnarray*}
 is the maximal solution to inequality (\ref{ZY2015:62}) out of the class of $[0,1]$-valued measurable functions on $[0,\infty)\times S.$ The statement immediately follows. (Note that there is a bounded nonnegative (measurable) solution to inequality (\ref{ZY2015:62}) if and only if there is a $[0,1]$-valued solution to it.) $\hfill\Box$\bigskip

The condition for explosion of a nonhomogeneous Markov pure jump process in a Borel state space in the above statement seems not stated in the literature. On the other hand, part (a) of Theorem \ref{ZY2015Theorem05} leads to a sharp drift condition for the nonexplosion, which is formulated as follows.
\begin{condition}\label{ZY2015Condition2}
There exists a monotone nondecreasing sequence $(\tilde{S}_n)\subseteq {\cal B}([0,\infty)\times S)$ and a $[0,\infty)$-valued measurable function $V$ on $[0,\infty)\times S$ such that the following hold.
\begin{itemize}
\item[(a)]  As $n\uparrow \infty,$ $\tilde{S}_n\uparrow [0,\infty)\times S$.
\item[(b)] For each $n=1,2,\dots,$
\begin{eqnarray}\label{ZY2015:57}
\sup_{x\in S_n,~t\ge 0} q_x(t)<\infty,
 \end{eqnarray}
 where $S_n$ denotes the projection of $\tilde{S}_n$ on $S$.
\item[(c)] As $n\uparrow \infty,$
\begin{eqnarray}\label{ZY2015:65}
\inf_{(t,x)\in ([0,\infty)\times S)\setminus \tilde{S}_n} V(t,x)\uparrow \infty.
\end{eqnarray}
\item[(d)] For some $\alpha \in (0,\infty),$
\begin{eqnarray}\label{ZY2015:61}
\int_0^\infty\int_S V(t+v,y) e^{-\alpha t -\int_0^t q_x(s+v)ds} \tilde{q}(dy|x,t+v)dt\le V(v,x),~\forall~x\in S,~v\in[0,\infty).
\end{eqnarray}

\end{itemize}
\end{condition}

\begin{theorem}\label{ZY2015Theorem03}
Consider the $Q$-function $q$ on $S$. Then the following assertions hold.
\begin{itemize}
\item[(a)] The nonhomogeneous Markov pure jump process $(X_t)$ with the transition function $P_q$ is nonexplosive if Condition \ref{ZY2015Condition2} is satisfied.

\item[(b)] Suppose $S$ is a locally compact separable metric space; $\sup_{x\in \Gamma,~t\ge 0}q_x(t)<\infty$ for each compact subset $\Gamma$ of $S$; and that the transition probability $p^\alpha$ for the discrete-time Markov chain $(\tilde{x}_n,\tilde{t}_n)$ is weakly Feller, i.e.,
    \begin{eqnarray*}
    (v,x)\in  [0,\infty)\times S\rightarrow \int_{S}\int_0^\infty f(t,y)p^\alpha(dt\times dy|v,x)
     \end{eqnarray*}
     is continuous for each bounded continuous function $f$ on $[0,\infty)\times S$.  If the nonhomogeneous Markov pure jump provess $(X_t)$ is nonexplosive, then Condition \ref{ZY2015Condition2} holds, where one can take $\tilde{S}_n=T_n\times S_n$ with $T_n\subseteq [0,\infty)$ and $S_n\subseteq S$ being compact sets such that
     \begin{eqnarray*}
     \sup_{(t,x)\in T_n\times S_n}V(t,x)<\infty.
     \end{eqnarray*}
\end{itemize}
\end{theorem}

\par\noindent\textit{Proof.} (a) Suppose Condition \ref{ZY2015Condition2} is satisfied. Let us extend the definition of the function $V$ to $([0,\infty)\times S)\bigcup (\{\infty\}\times \{\Delta\})$ by putting $V(\infty,\Delta)=0.$ Firstly, we observe that for each $n=0,1,\dots,$ the set $\tilde{S}_n\bigcup (\{\infty\}\times \{\Delta\})$ is petite for the Markov chain $(\tilde{t}_n,\tilde{x}_n)$ with the transition probability $p^\alpha.$ Indeed, for each $(v,x)\in \tilde{S}_n$ for some $n$,
\begin{eqnarray*}
&&p^\alpha(\Gamma_1\times \Gamma_2|v,x)\ge \delta_{(\infty,\Delta)}(\Gamma_1\times \Gamma_2) \left( 1- \int_0^\infty e^{-\alpha t -\int_0^t q_{x}(s+v)ds} q_{x }(t+v)dt\right)\\
&=&\delta_{(\infty,\Delta)}(\Gamma_1\times \Gamma_2)\int_0^\infty \alpha e^{-\alpha t} e^{-\int_0^t q_x(s+v)ds}dt\\
&\ge&\delta_{(\infty,\Delta)}(\Gamma_1\times \Gamma_2) \frac{\alpha}{\alpha+\sup_{x\in S_n,~t\ge 0} \{q_x(t)\}},
\end{eqnarray*}
c.f. (\ref{ZY2015:57}). Remember the notation used here that $S_n$ is the projection of $\tilde{S}_n$ on $S$.
The above inequality also trivially holds when
$(v,x)=(\infty,\Delta).$ By the validity of the previous inequality and (\ref{ZY2015:57}), the set $\tilde{S}_n\bigcup (\{\infty\}\times \{\Delta\})$ is petite for the chain $(\tilde{t}_n,\tilde{x}_n),$ see p.121 of Meyn and Tweedie \cite{MeynTweedie:1993}. Next, we note that the function $V$ is bounded off petite sets, i.e, for each $m=0,1,\dots,$ the set
$\{(v,x)\in ([0,\infty)\times S)\bigcup (\{\infty\}\times \{\Delta\}): V(v,x)\le m \}$ is petite. Indeed, by Condition \ref{ZY2015Condition2}, there exists some $n$ such that
\begin{eqnarray*}
\{(v,x)\in ([0,\infty)\times S)\bigcup (\{\infty\}\times \{\Delta\}): V(v,x)\le m \}\subseteq \tilde{S}_n\bigcup (\{\infty\}\times \{\Delta\}).
\end{eqnarray*}
Since the set of the right hand side is petite as observed earlier, so is the set on the left hand side. The chain $(\tilde{t}_n,\tilde{x}_n)$ is $\psi$-irreducible with the (maximal) irreducibility measure $\psi$ concentrated on $\{\infty\}\times \{\Delta\}$. Now for part (a), one can refer to Theorem 8.4.3 of \cite{MeynTweedie:1993}, where we put $C=\{\infty\}\times \{\Delta\}$, which is clearly petite.  (The validity of the proof of Theorem 8.4.3 in \cite{MeynTweedie:1993} does not require the state space of the chain to be locally compact.)

(b) By the local compactness assumption of the Borel space $S$, there exist a monotone nondecreasing sequence of open precompact sets $(T^o_n\times S^o_n)\subseteq {\cal B}([0,\infty)\times S)$ such that $T^o_n\times S^o_n\uparrow [0,\infty)\times S$. Now the statement is a consequence of the proof of Theorem 9.4.2 in \cite{MeynTweedie:1993}. (There, one can put $C_0=A_0= \{\infty\}\times \{\Delta\} $, $A_i=(T^o_i\times S^o_i)\bigcup(\{\infty\}\times \{\Delta\})$, and replace $A_i$ by its closure in (9.37) and (9.38).) It remains to take $\tilde{S}_n=T_n\times S_n$ as the closure of $T^o_n\times S^o_n.$ $\hfill\Box$\bigskip

The previous theorem shows that the drift condition for nonexplosion (Condition \ref{ZY2015Condition2}) is sharp. To the best of our knowledge, this condition has not been reported in the previous literature for nonhomogeneous Markov pure jump processes.

As an application of Theorem \ref{ZY2015Theorem03}, we shall apply it to obtaining directly a more applicable sufficient condition for nonexplosion, see Condition \ref{ZY2015Condition7} below, which is related to but at the same time more general than the one originally obtained by J. Zheng and X. Zheng \cite{Zheng:1987} and recently re-obtained by Chow and Khasminskii \cite{Chow:2011} for the process in a denumerable state space with a continuous $Q$-function. More comments on this are after Corollary \ref{ZY2015Corollary1} below.

 \begin{condition}\label{ZY2015Condition7}
There exist a monotone nondecreasing sequence $(S_n)\subseteq {\cal B}(S)$ and a nonnegative measurable function $V$ on $[0,\infty)\times S$ such that the following assertions hold.
\begin{itemize}
\item[(a)]  As $n\uparrow \infty,$ $S_n \uparrow  S$.
\item[(b)] For each $T>0$, as $n\uparrow \infty,$
\begin{eqnarray*}
\inf_{v\in[0,T],~x\in S\setminus S_n} V(v,x)\uparrow \infty.
\end{eqnarray*}
\item[(c)] For each $T>0$, and $n=1,2,\dots,$
\begin{eqnarray*}
\sup_{v\in[0,T],~x\in S_n}q_x(v)<\infty.
\end{eqnarray*}
\item[(d)] For each $T>0$, there is some constant $\alpha_T\in(0,\infty)$ such that
\begin{eqnarray*}
\int_0^{T-v}\int_S V(t+v,y) e^{-\alpha_T t -\int_0^t q_x(s+v)ds} \tilde{q}(dy|x,t+v)dt\le V(v,x),~\forall~v\in[0,T],~x\in S.
\end{eqnarray*}
\end{itemize}
 \end{condition}

\begin{corollary}\label{ZY2015Corollary1}
Consider a $Q$-function $q$ on $S$. The nonhomogeneous Markov pure jump process $(X_t)$ with the transition function $P_q$ is nonexplosive if Condition \ref{ZY2015Condition7} is satisfied.
\end{corollary}

\par\noindent\textit{Proof.} For each fixed $T\in[0,\infty),$ let us consider the following $Q$-function on $S$
\begin{eqnarray}\label{ZY2015:66}
{}_Tq(\Gamma|x,t):= q(\Gamma|x,t)I\{t\in[0,T]\},~\forall~\Gamma\in{\cal B}(S).
\end{eqnarray}
According to Feller's construction, c.f., (\ref{ZY2015:07}) and (\ref{ZY2015:01}), one can see that in particular,
\begin{eqnarray}\label{ZY2015:67}
P_q(0,x,t,\Gamma)=P_{{}_Tq}(0,x,t,\Gamma),~\forall~x\in S,~t\in[0,T],~\Gamma\in {\cal B}(S).
\end{eqnarray}
Note that under Condition \ref{ZY2015Condition7}, the corresponding version of Condition \ref{ZY2015Condition2} for the $Q$-function ${}_Tq$ is satisfied for each $T\in[0,\infty).$ Here, keeping in mind the definition (\ref{ZY2015:66}), only the validity of the corresponding version of (\ref{ZY2015:65}) needs special explanations, as follows. For each fixed $T\in[0,\infty),$  we can modify the definition of $V$ by putting
\begin{eqnarray*}
V_T(v,x)=I\{v>T\}V(0,x)+V(v,x)I\{v\in[0,T]\}
\end{eqnarray*}
for each $x\in S.$ In this way, we obtain for each $T\ge 0,$
\begin{eqnarray*}
\inf_{x\in S\setminus S_n,~v\in[0,\infty)}V_T(v,x)=\inf_{x\in S\setminus S_n,~v\in[0,T]}V(v,x)\rightarrow \infty,~\forall~x\in S,
\end{eqnarray*}
 as $n\rightarrow\infty$ under Condition \ref{ZY2015Condition7}. Consequently, for each $T>0,$ the corresponding version of Condition \ref{ZY2015Condition2} for the $Q$-function ${}_Tq$ is satisfied by the function $V_T$ and the sequence of sets $\tilde{S}_n=[0,\infty)\times S_n.$
By Theorem \ref{ZY2015Theorem03}, for each $T\ge 0$ and $x\in S,$ $P_{{}_Tq}(0,x,T,S)=1$, and thus
\begin{eqnarray*}
P_{q}(0,x,T,S)=1,~\forall~x\in S,~T\ge 0,
\end{eqnarray*}
by (\ref{ZY2015:67}). Now the statement follows from the above equality and Lemma \ref{ZY2015:Lemma07}. $\hfill\Box$\bigskip

The drift or test functions $V$ appearing in Conditions \ref{ZY2015Condition2} and  \ref{ZY2015Condition7}, though denoted with the same notation, are not necessarily the same. In our opinion, parts (b) and (c) of Condition \ref{ZY2015Condition7} are easier for practical verifications than parts (b) and (c) of Condition \ref{ZY2015Condition2}.

As mentioned earlier, recently, Chow and Khasminskii \cite{Chow:2011} also studied the nonexplosion of a nonhomogeneous Markov pure jump process but in a denumerable state space $S=\mathbb{Z}=\{0,\pm1,\pm2,\dots\}$; following a different approach, they obtained a sufficient condition for nonexplosion, which is formulated below in the setup of the present paper for the ease of reference of comparisons, c.f., Theorem 1 of \cite{Chow:2011}.

\begin{condition}[Zheng and Zheng, and Chow and Khasminskii]\label{ZY2015Condition6}
\par\noindent

\begin{itemize}
\item[(a)] $S=\mathbb{Z}$, and the function
\begin{eqnarray*}
t\in[0,\infty)\rightarrow q(\{j\}|i,t)
 \end{eqnarray*}
 is continuous for each $i,j\in S=\mathbb{Z}.$
 \item [(b)] There exists a nonnegative function $V$ on $[0,\infty)\times \mathbb{Z}$ such that the following assertions hold.
 \begin{itemize}
\item[(i)]$v\in [0,\infty)\rightarrow V(v,i)$ is continuously differentiable for each $i\in \mathbb{Z}$.
\item[(ii)] For each $T\ge 0$ and $i\in \mathbb{Z}$,
\begin{eqnarray*}
\sup_{0\le v\le T} \sum_{j\in \mathbb{Z}}\tilde{q}(\{j\}|i,v)|V(v,j)|<\infty.
\end{eqnarray*}
\item[(iii)] For each $T\ge 0$, as $n\uparrow \infty,$
\begin{eqnarray*}
\inf_{v\in[0,T],~|i|>n}V(v,i)\uparrow\infty.
\end{eqnarray*}
\item[(iv)] For each $T\ge 0$, there is some constant $\alpha_T\in(0,\infty)$ such that
\begin{eqnarray}\label{ZY2015:60}
\frac{\partial V(v,i)}{\partial v} + \sum_{j\in \mathbb{Z}} q(\{j\}|i,v)V(v,j)\le \alpha_T V(v,i),~\forall~v\in[0,T],~i\in \mathbb{Z}.
\end{eqnarray}
\end{itemize}
\end{itemize}
\end{condition}
We refer the reader to \cite{Chow:2011} for several examples, where the nonexplosion was verified using this drift condition. In Condition \ref{ZY2015Condition6}, the continuity of the $Q$-function on $S$ was imposed, so that the test function $V$ in the same condition is in the domain of the (infinitesimal) generator of the process $(t,X_t).$

It might be worthwhile to mention the following fact, which we learnt from Professors Mufa Chen (Beijing Normal University) and Floske Spieksma (Leiden University). Condition \ref{ZY2015Condition6} sufficient for the nonexplosion was, to the best our knowledge, initially obtained by J.Zheng and X.Zheng in the 1980s (also for the case of a demunerable state space and continous $Q$-function like in \cite{Chow:2011}). However, the authors announced the result without proof in \cite{Zheng:1987}, which, in spite of being in English, has not been widely accessible. The corresponding proof seems only available in the Ph.D thesis of J.Zheng \cite{Zheng:1993} under the supervision of Mufa Chen, which was written in Chinese, and not widely accessible neither.

The following observation shows that Condition \ref{ZY2015Condition6} is stronger than Condition \ref{ZY2015Condition2} and Condition \ref{ZY2015Condition7}.
\begin{proposition}\label{ZY2015PropositionS2}
If Condition \ref{ZY2015Condition6} is satisfied, then so are Condition \ref{ZY2015Condition2} and Condition \ref{ZY2015Condition7}.
 \end{proposition}

\par\noindent\textit{Proof.}
Suppose that Condition \ref{ZY2015Condition6} is satisfied. Since the concerned process is in a denumerable state space, and the $Q$-function is continuous, Theorem \ref{ZY2015Theorem03} applies, so that the validity of Condition \ref{ZY2015Condition2} follows from the sufficiency of Condition \ref{ZY2015Condition6} for nonexplosion, see Theorem 1 of \cite{Chow:2011}. To see that Condition \ref{ZY2015Condition7} is satisfied, merely note that for each $T\ge 0,$
\begin{eqnarray*}
&&\int_0^{T-v} \sum_{j\in \mathbb{Z}} V(t+v,j) e^{-\alpha_T t -\int_0^t q_i(s+v)ds} \tilde{q}(\{j\}|i,t+v)dt\\
&=&\int_0^{T-v}\sum_{j\in \mathbb{Z}} V(t+v,j) e^{-\alpha_T t -\int_0^t q_i(s+v)ds} q(\{j\}|i,t+v)dt\\
&&+\int_0^{T-v} V(t+v,i) e^{-\alpha_T t -\int_0^t q_i(s+v)ds} q_i(t+v)dt\\
&\le& \int_0^{T-v}  e^{-\alpha_T t -\int_0^t q_i(s+v)ds} (\alpha_T+q_i(t+v)) V(t+v,i) dt-\int_0^{T-v} e^{-\alpha_T t -\int_0^t q_i(s+v)ds}  \left.\frac{\partial V(v,i)}{\partial v}\right|_{t+v}dt\\
&\le& V(v,i),~\forall~v\in[0,T],~i\in \mathbb{Z},
\end{eqnarray*}
where the first inequality is by (\ref{ZY2015:60}). $\hfill\Box$ \bigskip

It is less transparent how to show directly that Condition \ref{ZY2015Condition6} implies Condition \ref{ZY2015Condition2} without using Theorem \ref{ZY2015Theorem03}.

It is clear that Conditions \ref{ZY2015Condition2} and \ref{ZY2015Condition7} are weaker than the following condition.
\begin{condition}\label{ZY2015Condition5}
There exists a monotone nondecreasing sequence $(S_n)\subseteq {\cal B}(S)$ and a $[0,\infty)$-valued measurable function $V$ on $S$ such that the following hold.
\begin{itemize}
\item[(a)]  As $n\uparrow \infty,$ $S_n \uparrow S$.
\item[(b)] For each $n=0,1,\dots,$
$
\sup_{x\in S_n,~v\ge 0} q_x(v)<\infty.
$
\item[(c)] As $n\uparrow \infty,$ $\inf_{x\in  S\setminus S_n} V(x)\uparrow \infty.$
\item[(d)] For some constant $\alpha\in(0,\infty),$
\begin{eqnarray*}
\int_S V(y)q(dy|x,v)\le \alpha V(x),~\forall~x\in S,~v\in[0,\infty).
\end{eqnarray*}
\end{itemize}
\end{condition}
In the context of continuous-time Markov decision processes, the version of the above continuous was imposed to guarantee the controlled process to be nonexplosive under each control policy in e.g., \cite{PiunovskiyZhang:2014}. If the $Q$-function $q(dy|x,s)$ on $S$ does not depend on time $s\ge 0,$ Condition \ref{ZY2015Condition5} reduces to Mufa Chen's condition, see \cite{Chen:1986,Chen:2004}. In p.112 of \cite{Chen:2004} and p.4 of \cite{Chen:2015}, intuitive explanations were provided as to why this condition is necessary for the nonexplosion in the homogeneous denumerable case, in particular, for the single-birth process, see Remark 3.20 of \cite{Chen:2004}. In the nonhomogeneous case, it is less clear when Condition \ref{ZY2015Condition5} is necessary for the nonexplosion.

To end this section, we say a few more words about the homogeneous case, i.e., when the $Q$-function $q(dy|x,s)$ on $S$ does not depend on time $s\ge 0,$ as follows. Then $(X_t)$ is a homogeneous Markov pure jump process with the transition function $P_q.$ Now instead of the discrete-time Markov chain $(\tilde{t}_n,\tilde{x}_n)$ with the transition probability $p_\alpha$ given by (\ref{ZY2015:53}), one can simply focus on the Markov chain $(\tilde{x}_n),$ whose transition probability can be obtained automatically from $p_\alpha$. The reasoning for the nonhomogeneous case in the previous discussions all applies. In particular, we see the following versions of Theorem \ref{ZY2015Theorem03}, Theorem \ref{ZY2015Theorem05} and Corollary \ref{ZY2015Corollary2} for the nonexplosion and nonexplosion of homogeneous Markov pure jump process $(X_t)$.
\begin{theorem}\label{ZY2015Theorem06}
Consider the homogeneous $Q$-function $q$ on $S$. Then the following assertions hold.
\begin{itemize}
\item[(a)] The homogeneous Markov pure jump process $(X_t)$ with the transition function $P_q$ is nonexplosive if Condition \ref{ZY2015Condition5} is satisfied, where the $Q$-function $q$ does not depend on time $v.$

\item[(b)] Suppose $S$ is a locally compact separable metric space, and the functions $x \in S\rightarrow q_x$ and $x\in S\rightarrow \int_S f(y)q(dy|x)$ are continuous for each bounded continuous function $f$ on $S$.
   If the homogeneous Markov pure jump process $(X_t)$ is nonexplosive, then Condition \ref{ZY2015Condition5} holds, where the sets $(S_n)$ can be taken as compact sets in $S$ such that
   \begin{eqnarray*}
   \sup_{x\in S_n}V(x)<\infty,~\forall~n=1,2,\dots.
   \end{eqnarray*}
\item[(c)] For each $\alpha\in (0,\infty)$, there is no nontrivial $[0,1]$-valued measurable function $U$ on $[0,\infty)\times S$ satisfying
 \begin{eqnarray*}
\alpha U(x)=\int_S U(y) q(dy|x),~\forall~x\in S.
 \end{eqnarray*}
\item[(d)] The process $(X_t)$ is explosive if and only if there exists a nontrivial bounded nonnegative measurable function $U$ on $S$ and some $\alpha>0$ such that
    \begin{eqnarray*}
    \alpha U(x)\le \int_S U(y) q(dy|x),~\forall~x\in S.
    \end{eqnarray*}
\end{itemize}
\end{theorem}
The first part of the above theorem is the same as Theorem 2.25 of \cite{Chen:2004}. Its origin is Theorem (16) in \cite{Chen:1986}, where a minor misprint was present, was obtained by Mufa Chen, using a different method. In the survey \cite{Chen:2015} and the monograph \cite{Chen:2005} (see p.167-168 therein) by Mufa Chen, that condition was demonstrated with various examples to be practically applicable. Part (c) of the above theorem coincides with Theorem 2.40 in \cite{Chen:2004}, where this necessary and sufficient condition for nonexplosion is called zero-exit. The origin of the zero-exit condition is in Theorem 7 of \cite{Reuter:1957} obtained by Reuter for the case of a homogeneous continuous-time Markov chain, for which, more discussions are given in \cite{Chen:2015,Spieksma:2015}.  Finally, part (d) of the previous statement was also presented in \cite{Chen:2004}, see Theorem 2.27 therein.

Part (b) of the above theorem was not reported elsewhere. On the other hand, also specialised to the homogeneous case but with a denumerable state space, the second part of Theorem \ref{ZY2015Theorem06}, asserting the necessity of Condition \ref{ZY2015Condition5} for the nonexplosion was established by Spieksma in \cite{Spieksma:2015}. In the present paper, we explored the idea in \cite{Spieksma:2015}. For the case of a homogeneous pure jump process in a general state space, the necessity of Condition \ref{ZY2015Condition5} was mentioned to be unclear in the recent survey \cite{Chen:2015}, see the proof of Theorem 4 therein, as well as the last paragraph of the text in p.223 of \cite{Chen:2015}. Theorem \ref{ZY2015Theorem06}(b) provides an affirmative answer when the state space is in addition locally compact, and the $Q$-function satisfies some weak Feller-type condition.

\section{Applicability of Dynkin's formula}\label{ZY2015:SecAppl}
In this section, we investigate when Dynkin's formula applies to a drift function, or more generally, to a measurable function bounded by a drift function.
\begin{definition}
Let $c\in[0,\infty)$ be fixed. A function $f$ on $S$ is called a $c$-drift function with respect to the $Q$-function $q$ on $S$ if
\begin{itemize}
\item[(a)]
it is $(0,\infty)$-valued and measurable on $S$, and
\item[(b)]  for each $x\in S,$
\begin{eqnarray}\label{ZY2015:37}
\int_{S}f(y)q(dy|x,s)\le c f(x)
\end{eqnarray}
for all $s\ge 0.$
\end{itemize}
\end{definition}

\begin{remark}
In (\ref{ZY2015:37}) of the previous definition, it is implicitly assumed that
\begin{eqnarray*}
\int_{S}f(y)\tilde{q}(dy|x,s)<\infty
\end{eqnarray*}
so that
the integral
\begin{eqnarray*}
\int_{S}f(y)q(dy|x,s):=\int_S f(y)\tilde{q}(dy|x,s)-q_x(s)f(x),
\end{eqnarray*}
i.e., the left hand side of (\ref{ZY2015:37}),
is well-defined and actually finite.
\end{remark}

Throughout the rest of this paper, when talking about a $c$-drift function, the constant $c$ is always nonnegative finite. In particular, each positive constant function is a $c$-drift function.

\begin{definition}
Let $f$ be a $c$-drift function with respect to the $Q$-function $q$ on $S$. A (real-valued) function $g$ on $S$ is said to be $f$-bounded if $g$ is measurable on $S$ and
\begin{eqnarray*}
||g||_f:=\sup_{x\in S}\frac{|g(x)|}{f(x)}<\infty.
\end{eqnarray*}
\end{definition}
The following fact is used frequently below without special reference.
 \begin{lemma}\label{ZY2015:Lemma08}
 Let $f$ be a $c$-drift function with respect to the $Q$-function $q$ on $S$. Then for each $f$-bounded function $g$ on $S$ and for each $x\in S$ and $s\ge 0$, the integral $\int_S g(y)q(dy|x,s)$ is well defined, and satisfies
 \begin{eqnarray*}\label{ZY2015:17}
\left|\int_S g(y)q(dy|x,s)\right|\le ||g||_f f(x)(c+2q_x(s)).
\end{eqnarray*}
\end{lemma}
\par\noindent\textit{Proof.} See the appendix. $\hfill\Box$\bigskip

\begin{remark}
For each $c$-drift function $f$ with respect to the $Q$-function $q$ on $S$, we automatically extend its definition to $S_\infty$ by putting \begin{eqnarray*}
f(x_\infty):=0.
\end{eqnarray*}
The similar extension is done to each $f$-bounded function $g$.
\end{remark}
According to the previous remark, in what follows, we write expressions such as $\int_S f(y)P_q(dy|0,x,t,dy)$ and $\mathbb{E}_x^q[f(X_t)]$ interchangeably.

\begin{lemma}\label{ZY2015:Lemma02}
For a $c$-drift function $f$ with respect to the $Q$-function $q$ on $S$,
\begin{eqnarray*}
\int_S f(y)P_q(dy|s,x,t,dy)\le e^{c(t-s)}f(x),~\forall~x\in S,~s,t\in [0,\infty), ~s\le t.
\end{eqnarray*}
Consequently,
\begin{eqnarray*}
\int_u^t \mathbb{E}_x^q\left[\left(\int_S f(y)q(dy|X_s,s)\right)^+\right]ds <\infty,~\forall~ x \in\ S,~u,t\in [0,\infty),~u\le t.
\end{eqnarray*}
\end{lemma}
\par\noindent\textit{Proof.} See the appendix. $\hfill\Box$\bigskip

The question of interest is when does the following Dynkin's formula apply to a $c$-drift function $f$ on $S$;
\begin{eqnarray}\label{ZY2015:04}
\mathbb{E}_x^q[f(X_t)]-f(x)=\int_0^t \mathbb{E}_x^q\left[\int_S f(y)q(dy|X_s,s)\right]ds,~\forall~ x \in S,~ t\in [0,\infty).
\end{eqnarray}
Note that by Lemma \ref{ZY2015:Lemma02}, for a $c$-drift function $f$, the right hand side is well defined, and in fact finite. More precisely, we have the next statement.

\begin{lemma}\label{ZY2015:Lemma01}
For a $c$-drift function $f$ with respect to the $Q$-function $q$ on $S$, relation (\ref{ZY2015:04}) holds if and only if
\begin{eqnarray}\label{ZY2015:03}
\mathbb{E}_x^q[f(X_t)]-\mathbb{E}_x^q[f(X_u)]=\int_u^t \mathbb{E}_x^q\left[\int_S f(y)q(dy|X_s,s)\right]ds,~\forall x \in\ S,~u,t\in [0,\infty),~u\le t.
\end{eqnarray}
Provided that either (\ref{ZY2015:04}) or (\ref{ZY2015:03}) holds, we have
\begin{eqnarray}\label{ZY2015:02}
\int_0^t \mathbb{E}^q_x\left[\left|\int_S f(y)q(dy|X_s,s)\right|\right]ds<\infty, ~\forall~x \in S,~ t\in [0,\infty).
\end{eqnarray}
\end{lemma}
\par\noindent\textit{Proof.} See the appendix. $\hfill\Box$\bigskip

The first assertion of the previous lemma is only needed in the proof of Lemma \ref{ZY2015:Lemma03} in the appendix.

\begin{theorem}\label{ZY2015Theorem02}
Let $f$ be a $c$-drift function with respect to the $Q$-function $q$ on $S$. Then
relation (\ref{ZY2015:04}) holds if and only if
\begin{eqnarray}\label{ZY2015:12}
P_{q^f}(0,x,t,S_\delta)=1,~\forall~x\in S,~t\ge 0,
\end{eqnarray}
i.e., the process $(X_t)$ with the transition function $P_{q^f}$ is nonexplosive (c.f., Lemma \ref{ZY2015:Lemma07}, and recall that the state $\delta$ is absorbing for the $Q$-function $q^f)$.
\end{theorem}

\par\noindent\textit{Proof.} Let $x\in S$ be arbitrarily fixed.
Then for each $t\ge 0,$
\begin{eqnarray}\label{ZY2015:10}
&&\int_S P_{q^f}(0,x,t,dy) q^f(S|y,t)= \int_S P_{q^f}(0,x,t,dy) \left( \frac{1}{f(y)} \int_S f(z)q(dz|y,t)-c\right)\nonumber\\
&=& \int_S \frac{e^{-ct}}{f(x)} f(y)P_q(0,x,t,dy) \left( \frac{1}{f(y)} \int_S f(z)q(dz|y,t)-c\right)\nonumber\\
&=&\frac{e^{-ct}}{f(x)} \left\{\int_S  \left(\int_S f(z)q(dz|y,t)\right)P_q(0,x,t,dy) - c\int_S f(y)P_q(0,x,t,dy)\right\},
\end{eqnarray}
where the first equality is by (\ref{ZY2015:08}), the second equality is by Lemma \ref{ZY2015:09}, and the last equality holds because the second integral inside the bracket therein is finite by Lemma \ref{ZY2015:Lemma02}.

Suppose (\ref{ZY2015:04}) holds.
By Lemma \ref{ZY2015:Lemma03},
\begin{eqnarray*}
&&\frac{e^{-ct}\mathbb{E}^q_x[f(X_t)]}{f(x)}=1+ \int_0^t \frac{e^{-c u}}{f(x)} \left( \mathbb{E}^q_x\left[ \int_S f(y)q(dy|X_u,u)\right]-c \mathbb{E}^q_x[f(X_u)]\right)du\\
&=&1+\int_0^t  \left(\int_S P_{q^f}(0,x,u,dy) q^f(S|y,u)\right)du,
\end{eqnarray*}
where the second equality is by (\ref{ZY2015:10}). Applying Lemma \ref{ZY2015:09} to the left hand side of the last equality, we see
\begin{eqnarray}\label{ZY2015:11}
P_{q^f}(0,x,t,S)=1+\int_0^t  \left(\int_S P_{q^f}(0,x,u,dy) q^f(S|y,u)\right)du.
\end{eqnarray}
On the other hand, by Proposition \ref{ZY2015:Prop03} applied to the (conservative and stable) $Q$-function $q^f$,
\begin{eqnarray*}
&&P_{q^f}(0,x,t,\{\delta\})=\int_0^t  \left(\int_S P_{q^f}(0,x,u,dy) q^f(\{\delta\}|y,u)\right)du\\
&=&\int_0^t  \left(\int_S P_{q^f}(0,x,u,dy) (-q^f(S|y,u))\right)du \in [0,1].
\end{eqnarray*}
Adding this equality to (\ref{ZY2015:11}), we see
\begin{eqnarray*}
P_{q^f}(0,x,t,S_\delta)=1,~\forall~x\in S, ~t\ge 0.
\end{eqnarray*}
This verifies (\ref{ZY2015:12}).

Suppose (\ref{ZY2015:12}) holds. Note that
\begin{eqnarray*}
\int_{S_\delta}(1+I\{y\in S\})q^f(dy|x,s)=q^f(S|x,s)\le 0\le c (I\{x\in S\}+1),~\forall~x\in S_\delta,~ s\ge 0,
\end{eqnarray*}
i.e., $y\rightarrow 1+I\{y\in S\}$ is a $c$-drift function with respect to the $Q$-function $q^f$ on $S_\delta.$
By (\ref{ZY2015:12}),
\begin{eqnarray*}
P_{q^f}(0,x,t,S_\delta)=1= \delta_x(S_\delta)+ \int_0^t \left(\int_S q^f(S_\delta|y,u)P_{q^f}(0,x,u,dy)\right)du,~\forall~t\ge 0,~x\in S,
\end{eqnarray*}
where the second equality holds because of the conservativeness of the $Q$-function $q^f.$
By Proposition \ref{ZY2015:Prop03} applied to $\Gamma=\{\delta\}$, the above equality leads to
\begin{eqnarray*}
P_{q^f}(0,x,t,S)= \delta_x(S)+ \int_0^t \left(\int_S q^f(S|y,u)P_{q^f}(0,x,u,dy)\right)du,~\forall~t\ge 0,~x\in S.
\end{eqnarray*}
This and (\ref{ZY2015:12}), which is assumed to be true, imply that (\ref{ZY2015:04}) holds for the $c$-drift function $y\rightarrow 1+I\{y\in S\}$ with respect to the $Q$-function $q^f$ on $S_\delta$, where $\mathbb{E}^q_x$ is replaced with $\mathbb{E}_x^{q^f}$, and $S$ is replaced by $S_\delta$. By Lemma \ref{ZY2015:Lemma03}, after simple calculations we see that for each $c>0,$
\begin{eqnarray*}
&&e^{c t} \mathbb{E}^{q^f}_x[I\{X_t\in S\}]=I\{x\in S\}+\int_0^t e^{c u} \left( \mathbb{E}_x^{q^f}\left[ \int_S I\{y\in S\}q^f(dy|X_u,u)\right]+c \mathbb{E}^{q^f}_x[I\{X_u\in S\}]\right)du,\\
&&~\forall~x\in S,~t\ge 0.
\end{eqnarray*}
Applying Lemma \ref{ZY2015:09} to the left hand side, the above equality reads
\begin{eqnarray*}
&& \frac{\mathbb{E}^q_x[f(X_t)]}{f(x)}=1+\int_0^t e^{c u} \left( \mathbb{E}_x^{q^f}\left[ q^f(S|X_u,u)\right]+c \mathbb{E}^{q^f}_x[I\{X_u\in S\}]\right)du,\\
&=&1+\int_0^t e^{c u} \left\{ \frac{e^{-cu}}{f(x)}\left( \mathbb{E}_x^{q}\left[ \int_S f(y)q(dy|X_u,u)\right]-c\int_S f(y)P_q(0,x,u,dy) \right)+c \mathbb{E}^{q^f}_x[I\{X_u\in S\}]\right\}du,\\
&=&1+\int_0^t e^{c u} \left\{ \frac{e^{-cu}}{f(x)}\left( \mathbb{E}_x^{q}\left[ \int_S f(y)q(dy|X_u,u)\right]-c\int_S f(y)P_q(0,x,u,dy) \right)\right.\\
&&\left.+c \frac{e^{-cu}}{f(x)}\int_S f(y)P_q(0,x,u,dy)\right\}du,\\
&=&1+ \int_0^t \frac{ \mathbb{E}_x^{q}\left[ \int_S f(y)q(dy|X_u,u)\right]}{f(x)}du, ~\forall~x\in S,~t\ge 0.
\end{eqnarray*}
where the second equality is by (\ref{ZY2015:10}), and the third equality is by Lemma \ref{ZY2015:09}. Now (\ref{ZY2015:04}) follows.
$\hfill\Box$
\bigskip

As an application of Theorem \ref{ZY2015Theorem02}, we consider the homogeneous case, and provide a condition for a $c$-drfit function $f$ to be in the domain of the extended generator. Here the term of extended generator is understood in the sense of \cite{MeynTweedie:1993III}.

\begin{corollary}
Consider a homogeneous $Q$-function $q$ on the Borel state space $S$, i.e., $q(dy|x,s)\equiv q(dy|x).$ A $c$-drift function $f$ is in the domain of the extended generator of the process $(X_t)$ with the transition function $P_q$ if there exists a monotone nondecreasing sequence $(S_n)\subseteq {\cal B}(S)$ and a $[0,\infty)$-valued measurable function $V$ on $S$ such that
Condition \ref{ZY2015Condition5}(a), (b) and (c) are satisfied, and
\begin{eqnarray*}
\int_S f(y)V(y)q(dy|x)\le \alpha f(x)V(x),~\forall~x\in S
\end{eqnarray*}
for some $\alpha>0.$ If in addition, $S$ is locally compact, and the kernel $\tilde{q}(dy|x)$ is weakly Feller, then the previous conditions are necessary.
\end{corollary}
\par\noindent\textit{Proof.}  Following the definition in \cite{MeynTweedie:1993III}, the $c$-drift function $f$ is in the domain of the extended generator of the process $(X_t)$ with the transition function $P_q$ if the relations (\ref{ZY2015:02}) and (\ref{ZY2015:04}) hold. Now the statement follows from Theorem \ref{ZY2015Theorem02} Lemma \ref{ZY2015:Lemma01} and Theorem \ref{ZY2015Theorem06}(a) and (b). $\hfill\Box$\bigskip

Next we generalize Theorem \ref{ZY2015Theorem02}. For each monotone nondecreasing sequence of subsets $(S_m)\subseteq {\cal B}(S)$, define the truncated $Q$-functions $q^{(m)}$ on $S$ by
\begin{eqnarray*}
q^m(\Gamma|x,s):=q(\Gamma|x,s)I\{x\in S_m\},~\forall~x\in S,~s\ge 0,~\Gamma\in {\cal B}(S)
\end{eqnarray*}
for each $m=0,1,\dots.$
Clearly, each state $x\in S\setminus S_m$, if exists, is absorbing with respect to the $Q$-function $q^{(m)}$. Without special reference, in what follows, by $(S_m)$ is always meant such a monotone nondecreasing sequence of subsets $(S_m)\subseteq {\cal B}(S)$ satisfying $S_m\uparrow S.$

\begin{theorem}\label{ZY2015Theorem01}
Consider the $Q$-functions $q$ and $(q^{(m)})$ on $S$. Let $x\in S,$ $s,t\ge 0$, $s\le t$ be arbitrarily fixed.
\par\noindent(a)
For each $\Gamma\in {\cal B}(S_j)$ for some $j=0,1,\dots,$ it holds that
\begin{eqnarray*}
\lim_{m\rightarrow \infty}P_{q^{(m)}}(s,x,t,\Gamma)= P_q(s,x,t,\Gamma).
\end{eqnarray*}
In fact, for each $M=0,1,\dots$ such that
$x\in S_M$ and $\Gamma \in {\cal B}(S_M)$,
\begin{eqnarray*}
P_{q^{(m)}}(s,x,t,\Gamma)\uparrow P_q(s,x,t,\Gamma)
\end{eqnarray*}
as $M\le m\uparrow \infty.$
\par\noindent(b) For each $\Gamma\in {\cal B}(S)$,
\begin{eqnarray*}
\liminf_{m\rightarrow \infty}P_{q^{(m)}}(s,x,t,\Gamma)\ge P_q(s,x,t,\Gamma).
\end{eqnarray*}
\end{theorem}

\par\noindent\textit{Proof.} (a) If $S_m=S$ for some $m=0,1,\dots,$ then the statement holds automatically. Below we consider $S_m\ne S$ for all $m=0,1,\dots.$ Firstly, note that for each $m=0,1,\dots,$
\begin{eqnarray}\label{ZY2015:15}
P_{q^{(m)}}(s,x,t,\Gamma)=\delta_x(\Gamma)=P_{q^{(m)}}^{(0)}(s,x,t,\Gamma),~\forall~x\in S\setminus S_m,~\Gamma\in {\cal B}(S),
\end{eqnarray}
since $\sup_{s\ge 0}\left\{q^{(m)}_x(s)\right\}=0$ for each $x\in S\setminus S_m;$ recall (\ref{ZY2015:07}).

Secondly, let us show that for each $m=0,1,\dots,$
\begin{eqnarray}\label{ZY2015:13}
\sum_{i=0}^n P_{q^{(m)}}^{(i)}(s,x,t,\Gamma)\le \sum_{i=0}^n P_{q^{(m+1)}}^{(i)}(s,x,t,\Gamma),~\forall~x\in S_m,~\Gamma\in {\cal B}(S_m),~s,t\ge 0,~s\le t,
\end{eqnarray}
by induction as follows. The case of $n=0$ is easily verified. Suppose (\ref{ZY2015:13}) holds for $n\le k$, and consider when $n=k+1.$ For each $x\in S_m,~\Gamma\in {\cal B}(S_m),~s,t\ge 0,~s\le t$,
\begin{eqnarray*}
&&\sum_{i=0}^{k+1} P_{q^{(m)}}^{(i)}(s,x,t,\Gamma) = P_{q^{(m)}}^{(0)}(s,x,t,\Gamma)+ \sum_{i=1}^{k+1} P_{q^{(m)}}^{(i)}(s,x,t,\Gamma)\\
&\le&P_{q^{(m+1)}}^{(0)}(s,x,t,\Gamma)+\sum_{i=0}^{k} P_{q^{(m)}}^{(i+1)}(s,x,t,\Gamma)\\
&=&P_{q^{(m+1)}}^{(0)}(s,x,t,\Gamma)+\sum_{i=0}^{k} \int_s^t e^{-\int_s^u q_x^{(m)}(v)dv} \int_S \tilde{q}^{(m)}(dy|x,u) P_{q^{(m)}}^{(i)}(u,y,t,\Gamma)\\
&=&P_{q^{(m+1)}}^{(0)}(s,x,t,\Gamma)+ \int_s^t e^{-\int_s^u q_x^{(m+1)}(v)dv} \left(\int_{S_m} \tilde{q}^{(m+1)}(dy|x,u) \sum_{i=0}^{k}P_{q^{(m)}}^{(i)}(u,y,t,\Gamma)\right.\\
&&\left.+\int_{S\setminus S_m} \tilde{q}^{(m+1)}(dy|x,u) \sum_{i=0}^{k}P_{q^{(m)}}^{(i)}(u,y,t,\Gamma)\right)\\
&\le &P_{q^{(m+1)}}^{(0)}(s,x,t,\Gamma)+\int_s^t e^{-\int_s^u q_x^{(m+1)}(v)dv} \int_S \tilde{q}^{(m+1)}(dy|x,u) \sum_{i=0}^{k} P_{q^{(m+1)}}^{(i)}(u,y,t,\Gamma)\\
&=&\sum_{i=0}^{k+1} P_{q^{(m+1)}}^{(i)}(s,x,t,\Gamma)
\end{eqnarray*}
where the first inequality is by the inductive supposition, the second equality is by (\ref{ZY2015:07}), the second inequality is by (\ref{ZY2015:15}), the definition of the $Q$-function $q^{(m)}$ and the induction supposition, and the last equality is by (\ref{ZY2015:07}). Thus, (\ref{ZY2015:13}) holds for the case of $n=k+1$, and by induction, for all $n=0,1,\dots.$

Thirdly, let us verify that
\begin{eqnarray}\label{ZY2015:16}
\lim_{m\rightarrow\infty} P_{q^{(m)}}^{(i)}(s,x,t,\Gamma)=P^{(i)}_{q}(s,x,t,\Gamma),~\forall~\Gamma\in {\cal B}(S),~x\in S,~s,t\ge0,~s\le t,
\end{eqnarray}
for each $i=0,1,\dots.$ Indeed, the case of $i=0$ is trivially true. Suppose (\ref{ZY2015:16}) holds for $i\le k$. Then
by (\ref{ZY2015:07}), the inductive supposition as well as the Lebesgue dominated convergence theorem, we see (\ref{ZY2015:16}) holds for $i=k+1,$ and thus for all $i=0,1,\dots.$

Now, let us fix some $x\in S,$ and $\Gamma\in {\cal B}(S_M)$ for some $M=0,1,\dots$. Without loss of generality, we assume $x\in S_M.$
Then
\begin{eqnarray*}
&&\lim_{m\rightarrow \infty}\lim_{n\rightarrow \infty} \sum_{i=0}^n P_{q^{(m)}}^{(i)}(s,x,t,\Gamma)=\lim_{n\rightarrow \infty}\lim_{m\rightarrow \infty} \sum_{i=0}^n P_{q^{(m)}}^{(i)}(s,x,t,\Gamma)\\
&=&\lim_{n\rightarrow \infty}\sum_{i=0}^n P_{q}^{(i)}(s,x,t,\Gamma)=P_q(s,x,t,\Gamma),
\end{eqnarray*}
where the first equality holds because $\sum_{i=0}^n P_{q^{(m)}}^{(i)}(s,x,t,\Gamma)$ increases in $n$ (resp. $m\ge M$) for each fixed $m$ (resp., $n$) by (\ref{ZY2015:13}), c.f. Theorem A1.6 of \cite{Bauerle:2011Book},
the second equality is by (\ref{ZY2015:16}), and the last equality is by (\ref{ZY2015:01}). The second assertion of this part of the statement holds because of the previous equality and (\ref{ZY2015:01}). The first assertion immediately follows from this.

(b) For each $x\in S,~s,~t\ge~0,~s\le t,$ and $\Gamma\in {\cal B}(S)$ and $m=0,1,\dots,$
\begin{eqnarray*}
&&P_{q^{(m)}}(s,x,t,\Gamma)=P_{q^{(m)}}\left(s,x,t,~\bigcup_{n=0}^\infty \left(\Gamma\bigcap (S_n\setminus S_{n-1})\right)\right)\\
&=& \sum_{n=0}^\infty P_{q^{(m)}}\left(s,x,t,~\Gamma\bigcap (S_n\setminus S_{n-1}) \right),
\end{eqnarray*}
where $S_0\setminus S_{-1}:=S_0;$ recall that $(S_n)$ is a monotone nondecreasing sequence.
By the Fatou lemma,
\begin{eqnarray*}
\liminf_{m\rightarrow \infty} P_{q^{(m)}}(s,x,t, \Gamma)\ge  \sum_{n=0}^\infty \liminf_{m\rightarrow\infty}P_{q^{(m)}}\left(s,x,t,~\Gamma\bigcap (S_n\setminus S_{n-1}) \right)=P_q(s,x,t,\Gamma),
\end{eqnarray*}
where the last equality is by part (a) of the statement.
$\hfill\Box$

\begin{lemma}\label{ZY2015:Lemma06}
Consider the $Q$-function $q$ on $S$. Let some $x\in S$ and $t\ge 0$ be arbitrarily fixed. Suppose $g$ is a measurable function on $S$ such that the following integral is finite:
\begin{eqnarray*}
\int_S g(y)q(dy|x,t)\in(-\infty,\infty),
\end{eqnarray*}
and
\begin{eqnarray}\label{ZY2015:18}
\mathbb{E}_x^q\left[\left|\int_S g(y)q(dy|X_t,t)\right|\right]<\infty.
\end{eqnarray}
Then
\begin{eqnarray*}
\lim_{m\rightarrow \infty} \mathbb{E}_x^{q^{(m)}}\left[ \int_S g(y)q^{(m)}(dy|X_t,t) \right]=\mathbb{E}_x^q\left[ \int_S g(y)q(dy|X_t,t) \right].
\end{eqnarray*}
\end{lemma}
\par\noindent\textit{Proof.} Let $m=0,1,\dots$ be fixed. Then
\begin{eqnarray}\label{ZY2015:33}
&&\mathbb{E}_x^{q^{(m)}}\left[\left|\int_S g(y)q^{(m)}(dy|X_t,t)\right|\right]=\mathbb{E}_x^{q^{(m)}}\left[\left|\int_S g(y)q(dy|X_t,t)I\{X_t\in S_m\}\right|\right]\nonumber\\
&=&\int_{S_m} P_{q^{(m)}}(0,x,t,dz)\left|\int_S g(y)q(dy|z,t)\right|
\le\int_{S_m} P_{q}(0,x,t,dz)\left|\int_S g(y)q(dy|z,t)\right|\nonumber\\
&\le&\int_{S} P_{q}(x,0,t,dz)\left|\int_S g(y)q(dy|z,t)\right|\\
&<&\infty,~\forall~x\in S_m\nonumber,
\end{eqnarray}
where the first inequality is by Theorem \ref{ZY2015Theorem01} (c.f., (\ref{ZY2015:13}) and Feller's construction), and the last inequality is by (\ref{ZY2015:18}). The above relation trivially holds for $x\in S\setminus S_m$, too. This fact and Theorem \ref{ZY2015Theorem01}(b) validate the application of Proposition \ref{ZY2015PropositionS1}; note that one can easily check from the definition of $q^{(m)}$ that
\begin{eqnarray*}
\lim_{m\rightarrow \infty} \int_S g(y)q^{(m)}(dy|x,t)=\int_S g(y)q(dy|x,t)
\end{eqnarray*}
for each $x\in S.$ It follows that the statement of this lemma holds.
$\hfill\Box$
\bigskip

\begin{lemma}\label{ZY2015Condition1}
There exists a monotone nondecreasing sequence of subsets $(S_m)_{m=0}^\infty\subseteq {\cal B}(S)$ such that
$
S_m\uparrow S,
$
and
$
\sup_{s\ge 0,~x\in S_m} \{q_x(s)\}<\infty,~\forall~m=0,1,\dots.
$
\end{lemma}
\par\noindent\textit{Proof.} See the proof of Lemma 1 of \cite{FeinbergShiryayev:2016}. $\hfill\Box$
\bigskip

\begin{theorem}\label{ZY2015Theorem07}
Let $f$ be a $c$-drift function for the $Q$-function $q$ on $S$. Then for
\begin{eqnarray}\label{ZY2015:68}
&&\mathbb{E}_x^{q}\left[g(X_t)\right]-g(x)=\int_0^t \mathbb{E}_x^{q}\left[\left(\int_S g(z)q(dz|X_s,s) \right) \right]ds\nonumber\\
&=& \mathbb{E}_x^{q}\left[\int_0^t \left(\int_S g(z)q(dz|X_s,s) \right) ds\right],~ \forall~x\in S,~t\ge0
\end{eqnarray}
\end{theorem}
to hold for each $f$-bounded function $g$ on $S$ satisfying
\begin{eqnarray}\label{ZY2015:29}
\int_0^t \mathbb{E}_x^{q}\left[\left|\int_S g(z)q(dz|X_t,t) \right| \right]dt<\infty,~ \forall~x\in S,~t\ge0,
\end{eqnarray}
it is necessary and sufficient that one of the following equivalent (by Theorem \ref{ZY2015Theorem02}) conditions hold:
\begin{itemize}
\item[(a)] The relation (\ref{ZY2015:04}) holds.
\item[(b)] The process $(X_t)$ with the transition function $P_{q^f}$ is nonexplosive.
\end{itemize}

\par\noindent\textit{Proof.} Throughout this proof, let $x\in S$ and $t\ge 0$ be arbitrarily fixed. Consider the sequence of sets $(S_m)$ in Lemma \ref{ZY2015Condition1} and the associated $Q$-functions $(q_m).$ Then, one can check that for each $m=0,1,\dots,$ the function $f$ is a $c$-drift function on $S$ for the $Q$-function $q_m$. As a result, for each $m=0,1,\dots,$
\begin{eqnarray}\label{ZY2015:36}
\mathbb{E}_x^{q^{(m)}}\left[g(X_t)\right]-g(x)=\int_0^t \mathbb{E}^{q^{(m)}}_x\left[\left(\int_S g(z)q^{(m)}(dz|X_s,s) \right) \right]ds
\end{eqnarray}
for each $f$-bounded function $g$ by Lemma \ref{ZY2015:Lemma05}; recall that $\sup_{x\in S,~s\ge 0}\{q^{(m)}_x(s)\}<\infty.$

The necessity part of the statement is clear because of Theorem \ref{ZY2015Theorem02} (see also Lemma \ref{ZY2015:Lemma01}). Therefore, we only prove the sufficiency part as follows. Consider the function $g$ that is $f$-bounded and satisfies (\ref{ZY2015:29}).
By (\ref{ZY2015:29}), for almost all $s\in[0,t]$ (with respect to the Lebesgue measure),
\begin{eqnarray*}
\mathbb{E}_x^{q}\left[\left|\int_S g(z)q(dz|X_s,s) \right| \right]<\infty.
\end{eqnarray*}
Lemma \ref{ZY2015:Lemma06} implies that
\begin{eqnarray*}
\lim_{m\rightarrow \infty} \mathbb{E}_x^{q^{(m)}}\left[ \int_S g(y)q^{(m)}(dy|X_s,s) \right]=\mathbb{E}_x^q\left[ \int_S g(y)q(dy|X_s,s) \right]
\end{eqnarray*}
for almost all $s\in[0,t].$ By (\ref{ZY2015:33}) with $t$ being replaced by $s\in[0,t]$, which holds for all $x\in S$ (as mentioned immediately after it), and (\ref{ZY2015:29}), one can apply the Lebesgue dominated convergence theorem for that
\begin{eqnarray}\label{ZY2015:35}
\lim_{m\rightarrow \infty}\int_0^t \mathbb{E}_x^{q^{(m)}}\left[\left(\int_S g(z)q^{(m)}(dz|X_s,s) \right) \right]ds= \int_0^t \mathbb{E}_x^{q}\left[\left(\int_S g(z)q(dz|X_s,s) \right) \right]ds.
\end{eqnarray}
Note that the above relation holds if $g$ is replaced by $f$, because (\ref{ZY2015:04}) is supposed to hold, and Lemma \ref{ZY2015:Lemma01}, see (\ref{ZY2015:02}) therein.

Next, let us show that
\begin{eqnarray}\label{ZY2015:32}
\lim_{m\rightarrow \infty}\mathbb{E}_x^{q^{(m)}}[g(X_t)]=\mathbb{E}^{q}_x[g(X_t)]
\end{eqnarray}
as follows. By Lemma \ref{ZY2015:09} applied to the $Q$-function ${q^{(m)}}^f$,
\begin{eqnarray}\label{ZY2015:30}
\frac{e^{-ct}}{f(x)}\int_{\Gamma} f(y)P_{q^{(m)}}(0,x,t,dy)=P_{{q^{(m)}}^f}(0,x,t,\Gamma),~\forall~\Gamma\in {\cal B}(S).
\end{eqnarray}
On the other hand,
\begin{eqnarray*}
&&\limsup_{m\rightarrow \infty} \int_\Gamma f(y) P_{q^{(m)}}(0,x,t,dy)=\limsup_{m\rightarrow \infty} \int_S f(y) I\{y\in \Gamma\} P_{q^{(m)}}(0,x,t,dy)\\
&\le& \limsup_{m\rightarrow \infty} \int_S f(y) P_{q^{(m)}}(0,x,t,dy)\\
&=&f(x)+\limsup_{m\rightarrow \infty}\int_0^t \int_S \left(\int_{S}f(z)q^{(m)}(dz|y,s) \right) P_{q^{(m)}}(0,x,s,dy)ds \\
&=&f(x)+\int_0^t \int_S \left(\int_{S}f(z)q(dz|y,s) \right) P_{q}(0,x,s,dy)ds\\
&=&\mathbb{E}^{q}_x[f(X_t)]\le e^{ct}(x)<\infty,~\forall~\Gamma\in {\cal B}(S),
\end{eqnarray*}
where the second equality holds because (\ref{ZY2015:36}) holds for $g=f$, and the third equality holds due to that (\ref{ZY2015:35}) is valid for $g=f$, see the remark immediately after (\ref{ZY2015:35}), the fourth equality is by validity of (\ref{ZY2015:04}) for the function $f$, and the second but the last inequality is by Lemma \ref{ZY2015:Lemma02}. Therefore, keeping in mind Theorem \ref{ZY2015Theorem01}(b) applied to the $Q$-function $q^{(m)}$,
one can refer to Proposition \ref{ZY2015PropositionS1} for that
\begin{eqnarray*}
\lim_{m\rightarrow \infty}  \int_\Gamma f(y) P_{q^{(m)}}(0,x,t,dy)=\int_\Gamma f(y)P_q(0,x,t,dy),~\forall~\Gamma\in {\cal B}(S).
\end{eqnarray*}
The above equality and (\ref{ZY2015:30}) imply
\begin{eqnarray}\label{ZY2015:31}
&&\lim_{m\rightarrow \infty}P_{{q^{(m)}}^f}(0,x,t,\Gamma)=\frac{e^{-ct}}{f(x)}\int_{\Gamma} f(y)P_{q}(0,x,t,dy)=P_{q^f}(0,x,t,\Gamma),~\forall~\Gamma\in {\cal B}(S),
\end{eqnarray}
where the second equality is by
Lemma \ref{ZY2015:09}. It follows from the previous setwise convergence that
\begin{eqnarray*}
&&\lim_{m\rightarrow \infty}\frac{e^{-ct}}{f(x)}\int_{S} g(y) P_{q^{(m)}}(0,x,t,dy)=\lim_{m\rightarrow \infty}\frac{e^{-ct}}{f(x)}\int_{S}\frac{g(y)}{f(y)} f(y)P_{q^{(m)}}(0,x,t,dy)\\
&=&\lim_{m\rightarrow \infty}\int_S P_{{q^{(m)}}^f}(0,x,t,dz)\frac{g(y)}{f(y)}=\frac{e^{-ct}}{f(x)}\int_{S} \frac{g(y)}{f(y)}f(y)P_{q}(0,x,t,dy) \\
&=&\frac{e^{-ct}}{f(x)}\int_{S}  g(y) P_{q}(0,x,t,dy),
\end{eqnarray*}
where the second equality is by (\ref{ZY2015:30}), and for the third equality we recall (\ref{ZY2015:31}) and that the measurable function $y\in S\rightarrow \frac{g(y)}{f(y)}$ is bounded on $S$. Consequently, (\ref{ZY2015:32}) holds.

Now the statement follows from (\ref{ZY2015:35}), (\ref{ZY2015:32}) and (\ref{ZY2015:36}). $\hfill\Box$\bigskip

In the denumerable and homogeneous case, the previous theorem was obtained in \cite{Spieksma:2015}. Here we applied the idea in \cite{Spieksma:2012} with some additional preparations; e.g., now we had to justify Theorem \ref{ZY2015Theorem01}, whereas for the homogeneous denumerable case Proposition 2.14 of Anderson \cite{Anderson:1991} was directly referred to in \cite{Spieksma:2015}. In that connection, we mention that although the relevant part to \cite{Spieksma:2012} in Proposition 2.14 of \cite{Anderson:1991} is correct, the general conclusion of Proposition 2.14 of \cite{Anderson:1991} does not hold; see \cite{Chen:1996}.

\section{Conclusion}\label{ZY2015ConcSec}

In this paper, we obtained conditions for nonexplosion and explosion for nonhomogeneous Markov pure jump processes in Borel state spaces. The conditions are sharp; e.g., the one for nonexplosion is necessary if the state space is locally compact and the $Q$-function satisfies weak Feller-type conditions. For the homogeneous case, some known conditions are naturally retrieved. The argument here extends the one in \cite{Spieksma:2015}, which deals with the (homogeneous) continuous-time Markov chains, and is different from those in \cite{Chow:2011,Zheng:1993} for the nonhomogeneous denumerable state space case, and \cite{Chen:1986,Chen:2004} for the homogeneous general state space case. In particular, this argument easily shows the sharpness of the obtained condition, which partially addresses the remark made by Mufa Chen in the recent survey \cite{Chen:2015}, see the proof of Theorem 4 as well as the last paragraph of the text in p.223 of \cite{Chen:2015}. Extending the argument in \cite{Spieksma:2012}, the relations between nonexplosion and the applicability of Dynkin's formula to a drift function were demonstrated, too. This type of result was used in \cite{Blok:2015} to study a controlled continuous-time Markov chain with stationary policies.

\appendix
\section{Appendix}
In this appendix, we present some auxiliary results and the proofs of some lemmas stated in the previous sections.

\begin{proposition}\label{ZY2015PropositionS1}
Let $\mu_n$ and $\mu$ be $[0,\infty]$-valued measures on the Borel space $X,$ and $f_n$ (resp., $f$) be $(-\infty,\infty)$-valued $\mu_n$-integrable (resp., measurable) functions thereon, and $g_n$ be nonnegative measurable functions on $X$.  Suppose $\liminf_{n\rightarrow \infty}\mu_n(\Gamma)\ge \mu(\Gamma)$ for each $\Gamma\in {\cal B}(X)$, and $f_n(x)\rightarrow f(x)$ for each $x\in X.$ If $|f_n(x)|\le g_n(x)$ for each $x\in X$ and $n\ge 1,$ and
\begin{eqnarray*}
\limsup_{n\rightarrow \infty} \int_X g_n(y)\mu_n(dy)\le \int_X \liminf_{n\rightarrow\infty}g_n(y)\mu(dy)<\infty,
\end{eqnarray*}
then $f$ is $\mu$-integrable, and
\begin{eqnarray*}
&&\lim_{n\rightarrow\infty} \int_X f_n(y)\mu_n(dy)=\int_X f(y)\mu(dy)\in (-\infty,\infty).
\end{eqnarray*}
\end{proposition}
\par\noindent\textit{Proof.} See Theorem 2.4 of \cite{Serfozo:1982} $\hfill\Box$\bigskip

\par\noindent\textit{Proof of Lemma \ref{ZY2015:Lemma07}.} Let us firstly justify that (b), (c) and (d) are equivalent. Clearly, (d) implies (b), which in turn implies (c). That (c) implies (b) is true because clearly $t\rightarrow P(0,x,t,S)$ is monotone nonincreasing in $t\ge 0$. Now suppose (b) is true. By the Chapman-Kolmogorov equation, we see that for each $x\in S$ and $0\le s\le t$,
 \begin{eqnarray}\label{ZY2015:38}
 P_q(s,y,t,S)=1
 \end{eqnarray}
 for almost all $y$ with respect to $P_q(0,x,s,dy).$
Now suppose for contradiction that (d) does not hold, so that there exist some $0\le s<t$ and $x\in S$ such that
 \begin{eqnarray}\label{ZY2015:39}
 P_q(s,x,t,S)<1.
 \end{eqnarray}
Then \begin{eqnarray*}
&&1=P_q(0,x,t,S)=P_q(0,x,s,\{x\})P(s,x,t,S)+\int_{S\setminus \{x\}} P_q(0,x,s,dy)P_q(s,y,t,S)\\
&=&P_q(0,x,s,\{x\})P(s,x,t,S)+ P_q(0,x,s,S\setminus \{x\})\\
&<&P_q(0,x,s,\{x\})+ P_q(0,x,s,S\setminus \{x\})=1,
\end{eqnarray*}
where the first equality and the last equality are by assumption that (b) holds, the second equality is by the Chapman-Kolmogorov equation, the third equality is by (\ref{ZY2015:38}), and the inequality is by (\ref{ZY2015:39}). This is a desired contradiction, and thus assertions (b), (c) and (d) are equivalent.

It is clear that (a) implies (b), c.f. Proposition \ref{ZY2015:Prop01} and (\ref{ZY2015:50}). That (b) implies (a) follows from the relation
$
\{t_\infty=\infty\}=\bigcap_{n=1}^\infty \{X_n\in S\},
$
c.f. (\ref{ZY2015:50}). The proof is completed. $\hfill\Box$\bigskip

\par\noindent\textit{Proof of Lemma \ref{ZY2015:Lemma08}.}
It is clear that $\int_S g(y)q(dy|x,s)$ is well defined. For the second part of the statement, note that for each $x\in S$ and $s\ge 0$,
\begin{eqnarray*}
&&\int_S g(y)\tilde{q}(dy|x,s)-g(x)q_x(s)\ge -||g||_f \int_S f(y)\tilde{q}(dy|x,s)-||g||_f f(x)q_x(s)\\
&=&-||g||_f \int_S f(y)q(dy|x,s)- 2||g||_ff(x)q_x(s) \ge -||g||_f cf(x)- 2||g||_ff(x)q_x(s)\\
&=& -||g||_f f(x)(c+2q_x(s)),
\end{eqnarray*}
where the second inequality is by the fact that $f$ is a $c$-drift function. Thus,
\begin{eqnarray*}
\int_S g(y)q(dy|x,s)\ge  -||g||_f f(x)(c+2q_x(s)),
\end{eqnarray*}
the integral on the left hand side being well defined. Similarly, one can check that
\begin{eqnarray*}
\int_S g(y)q(dy|x,s)\le ||g||_f f(x)(c+2q_x(s)),
\end{eqnarray*}
as required. $\hfill\Box$\bigskip

\par\noindent\textit{Proof of Lemma \ref{ZY2015:Lemma02}.} The first assertion follows from the proof of the relevant assertion in Theorem 3.1 of \cite{Guo:2007}. The second assertion follows from the first assertion because for each $x\in S$ such that  $\int_{S}f(y)q(dy|x,s)\ge 0,$
\begin{eqnarray*}
\left(\int_S f(y)q(dy|x,s)\right)^+=\int_{S}f(y)q(dy|x,s) \le cf(x),~ \forall~s\ge 0,
\end{eqnarray*}
whereas for each  $x\in S$ such that  $\int_{S}f(y)q(dy|x,s)< 0,$
\begin{eqnarray*}
\left(\int_S f(y)q(dy|x,s)\right)^+=0 \le cf(x),~ \forall~s\ge 0.
\end{eqnarray*}
The statement is proved.
$\hfill\Box$

\begin{lemma}\label{ZY2015:Lemma05}
Let $f$ be a $c$-drift function with respect to the $Q$-function $q$ on $S$. Suppose
\begin{eqnarray*}
\sup_{x\in S,~s\ge 0}\{q_x(s)\}<\infty.
\end{eqnarray*}
Then for each $x\in S,$ $0\le s\le t$, and $f$-bounded (measurable) function $g$ on $S$,
\begin{eqnarray*}
\int_S g(z)P_q(s,x,t,dz)=g(x)+\int_s^t \left(\int_Sg(z)  q(dz|y,u)\right)P_q(s,x,u,dy)du.
\end{eqnarray*}
In particular,
\begin{eqnarray}\label{ZY2015:23}
P_q(s,x,t,S)=1,~\forall~x\in S, ~s,t\in [0,\infty),~s\le t.
\end{eqnarray}
\end{lemma}
\par\noindent\textit{Proof.} The first assertion immediately follows from Proposition \ref{ZY2015:Prop03} and Lemma \ref{ZY2015:Lemma02}. The last assertion follows because $f\equiv 1$ is a $c$-drift function, and one can put $g\equiv 1$ in the first assertion. $\hfill\Box$\bigskip

\par\noindent\textit{Proof of Lemma \ref{ZY2015:Lemma01}.} We only need justify that (\ref{ZY2015:04}) implies (\ref{ZY2015:03}) as follows. Let some $x \in S$ and $t\in [0,\infty)$ be arbitrarily fixed. Since (\ref{ZY2015:04}) holds by assumption,
\begin{eqnarray*}
\int_0^t \mathbb{E}^q_x\left[\left(\int_S f(y)q(dy|X_s,s)\right)^-\right]ds=\int_0^t \mathbb{E}^q_x\left[\left(\int_S f(y)q(dy|X_s,s)\right)^+\right]ds
-\mathbb{E}^q_x[f(X_t)]+f(x)<\infty
\end{eqnarray*}
by Lemma \ref{ZY2015:Lemma02}. Consequently, (\ref{ZY2015:02}) holds.
According to the Fubini theorem,
\begin{eqnarray*}
\int_0^t \mathbb{E}^q_x\left[\int_S f(y)q(dy|X_s,s)\right]ds= \mathbb{E}^q_x\left[\int_0^t \left(\int_S f(y)q(dy|X_s,s)\right)ds\right]\in(-\infty,\infty).
\end{eqnarray*}
Now for each $0\le u\le t,$ (\ref{ZY2015:04}) implies
\begin{eqnarray*}
&&\mathbb{E}_x^q[f(X_t)]-f(x)= \mathbb{E}_x^q\left[\int_0^t\int_S f(y)q(dy|X_s,s)ds\right]\\
&=&\mathbb{E}_x^q\left[\int_0^u\int_S f(y)q(dy|X_s,s)ds\right]+\mathbb{E}_x^q\left[\int_u^t\int_S f(y)q(dy|X_s,s)ds\right]\\
&=& \mathbb{E}_x^q[f(X_u)]-f(x)+\int_u^t\mathbb{E}_x^q\left[\int_S f(y)q(dy|X_s,s)\right]ds.
\end{eqnarray*}
All the expressions are finite. Now (\ref{ZY2015:03}) follows. $\hfill\Box$

\begin{lemma}\label{ZY2015:Lemma03}
For a $c$-drift function $f$ with respect to the $Q$-function on $S$ that satisfies relation (\ref{ZY2015:04}), it holds that for each $d\in (-\infty,\infty)$, $x\in S$ and $t\ge 0$,
\begin{eqnarray*}
e^{d t} \mathbb{E}^q_x[f(X_t)]=f(x)+\int_0^t e^{d u} \left( \mathbb{E}^q_x\left[ \int_S f(y)q(dy|X_u,u)\right]+d \mathbb{E}^q_x[f(X_u)]\right)du.
\end{eqnarray*}
\end{lemma}
\par\noindent\textit{Proof.} Keeping in mind Lemmas \ref{ZY2015:Lemma02} and \ref{ZY2015:Lemma01} (and especially (\ref{ZY2015:02})), direct calculations give
\begin{eqnarray*}
&&\int_0^t  e^{ds}\mathbb{E}_x^q\left[\int_S f(y)q(dy|X_s,s)\right]ds\\
&=&\int_0^t \int_0^s \left(d e^{du}du\right) \mathbb{E}^q_x\left[\int_S f(y)q(dy|X_s,s)\right]ds +\mathbb{E}^q_x[f(X_t)]-f(x)\\
&=& \int_0^t d e^{du} \left(\int_u^t \mathbb{E}^q_x\left[\int_S f(y)q(dy|X_s,s)\right]ds\right)du+\mathbb{E}^q_x[f(X_t)]-f(x)\\
&=& e^{dt} \mathbb{E}^q_x[f(X_t)](e^{dt}-1)- \int_0^t de^{du} \mathbb{E}^q_x[f(X_u)]du+\mathbb{E}^q_x[f(X_t)]-f(x),
\end{eqnarray*}
where the first and the third equalities are by (\ref{ZY2015:03}). All the expressions are finite. The statement now follows. $\hfill\Box$\bigskip

\begin{lemma}\label{ZY2015:09}
Consider a $c$-drift function $f$ with respect to the $Q$-function $q$ on $S.$ For each $x\in S,$ $s,t\ge 0$, $s\le t$ and $\Gamma\in {\cal B}(S)$, the following relation holds;
\begin{eqnarray*}
P_{q^f}(s,x,t,\Gamma)=\frac{e^{-c(t-s)}}{f(x)}\int_\Gamma f(y) P_q(s,x,t,dy).
\end{eqnarray*}
\end{lemma}
\par\noindent\textit{Proof.} By Feller's construction, c.f., (\ref{ZY2015:07}) and (\ref{ZY2015:01}), it suffices to show
\begin{eqnarray}\label{ZY2015:05}
\sum_{m=0}^n P_{q^f}^{(m)}(s,x,t,\Gamma)=\frac{e^{-c(t-s)}}{f(x)}\int_\Gamma f(y) \left(\sum_{m=0}^n P^{(m)}_q(s,x,t,dy)\right),~x\in S,~s,t\ge 0,~ s\le t,~\Gamma\in {\cal B}(S),
\end{eqnarray}
as follows. Let $x\in S,~s,t\ge 0,~ s\le t,~\Gamma\in {\cal B}(S)$ be arbitrarily fixed. Consider $n=0.$ Then
\begin{eqnarray}\label{ZY2015:06}
&&\frac{e^{-c(t-s)}}{f(x)}\int_\Gamma f(y)   P^{(0)}_q(s,x,t,dy)=\frac{e^{-c(t-s)}}{f(x)}\int_\Gamma f(y)  e^{-\int_s^t q_x(v)dv} \delta_x(dy) =\delta_x(\Gamma) e^{-\int_s^t q_x(v)dv-c(t-s)}\nonumber\\
&=&P_{q^f}^{(0)}(s,x,t,\Gamma).
\end{eqnarray}
Now assume that (\ref{ZY2015:05}) holds for $n=k$.
Then
\begin{eqnarray*}
&&\frac{e^{-c(t-s)}}{f(x)}\int_\Gamma f(y) \left(\sum_{m=0}^{k+1} P^{(m)}_q(s,x,t,dy)\right)\\
&=&P_{q^f}^{(0)}(s,x,t,\Gamma)+\frac{e^{-c(t-s)}}{f(x)}\int_\Gamma f(y) \left(\sum_{m=1}^{k+1} P^{(m)}_q(s,x,t,dy)\right)\\
&=&P_{q^f}^{(0)}(s,x,t,\Gamma)+\frac{e^{-c(t-s)}}{f(x)}\int_\Gamma f(y) \left(\sum_{m=0}^{k} P^{(m+1)}_q(s,x,t,dy)\right)\\
&=&P_{q^f}^{(0)}(s,x,t,\Gamma)+\frac{e^{-c(t-s)}}{f(x)}\int_\Gamma f(y) \left(\sum_{m=0}^{k}\int_s^t e^{-\int_s^u q_x(v)dv} \left(\int_S \tilde{q}(dz|x,u) P^{(m)}_q(u,z,t,dy)\right)du\right)\\
&=&P_{q^f}^{(0)}(s,x,t,\Gamma)+\int_{s}^t e^{-\int_s^u (q_x(v)+c)dv} \int_S \frac{f(z)}{f(x)} \tilde{q}(dz|x,u) \left(\int_\Gamma \frac{f(y)}{f(z)} e^{-c(t-u)}\sum_{m=0}^k P^{(m)}_q(u,z,t,dy) \right)du\\
&=&P_{q^f}^{(0)}(s,x,t,\Gamma)+\sum_{m=0}^k P_{q^f}^{(m+1)}(s,x,t,\Gamma),
\end{eqnarray*}
where the first equality is by (\ref{ZY2015:06}), the third equality is by (\ref{ZY2015:07}), and the last equality is by (\ref{ZY2015:07}), the inductive supposition and (\ref{ZY2015:08}). Now (\ref{ZY2015:05}) holds for $n=k+1$, and thus for all $n=0,1,\dots$ by induction. $\hfill\Box$

\par\noindent\textbf{Acknowledgement.} I would like to thank Professor Mufa Chen (Beijing Normal University) for providing the scan copy of the relevant pages in \cite{Zheng:1993} and the paper \cite{Zheng:1987}. I also thank the referees for their helpful comments and remarks. This work was carried out with a financial grant from the Research Fund for Coal and Steel of the European Commission, within the INDUSE-2-SAFETY project (Grant No. RFSR-CT-2014-00025).

\end{document}